# CONDITIONED SQUARE FUNCTIONS FOR NONCOMMUTATIVE MARTINGALES[1]

By Narcisse Randrianantoanina

*Miami University*

We prove a weak-type $(1,1)$ inequality involving conditioned versions of square functions for martingales in noncommutative $L^p$-spaces associated with finite von Neumann algebras. As application, we determine the optimal orders for the best constants in the noncommutative Burkholder/Rosenthal inequalities from [*Ann. Probab.* **31** (2003) 948–995]. We also discuss BMO-norms of sums of noncommuting order-independent operators.

**1. Introduction.** The role played by martingales in the development of classical probability and analysis is well known as evidenced in the books [5, 7, 8, 21]. In recent years, many classical inequalities from classical martingale theory have been reformulated to include noncommutative martingales. Several articles on the subject of noncommutative martingales have appeared in the literature recently. We refer the reader to a recent survey by Xu [35] for an up-to-date exposition of this topic.

In this paper, we continue this line of research by studying conditioned square functions of noncommutative martingales. Recall that conditioned square function inequalities evolved from a classical result of Rosenthal [32] on $p$-moment of sums of independent mean-zero random variables back in the 1970s which was later generalized by Burkholder for the context of martingales [2] as follows: let $2 \leq p < \infty$ and let $(\mathcal{F}_n)$ be a filtration on a probability space $(\Omega, \mathcal{F}, \mathbb{P})$. Given $f \in L^p$, the conditional expectations $(\mathcal{E}_n)$ and the martingale difference sequence are given by

$$\mathcal{E}_n(x) = E(f|\mathcal{F}_n) \quad \text{and} \quad d_n = \mathcal{E}_n(f) - \mathcal{E}_{n-1}(f).$$

Received October 2005; revised April 2006.
[1]Supported in part by NSF Grant DMS-04-56781.

Then the following equivalence holds:

$$\left(\sum_{n\geq 1} \|d_n\|_p^p\right)^{1/p} + \left\|\left(\sum_{n\geq 1}\mathcal{E}_{n-1}(|d_n|^2)\right)^{1/2}\right\|_p \sim_{c_p} \|f\|_p,$$

where $A \sim_c B$ means $c^{-1}A \leq B \leq cA$. This equivalence was inspired by Rosenthal's inequality for sums of independent mean-zero random variables $(f_n)$ as he obtained the above equivalence for $d_n = f_n$ and $\mathcal{E}_{n-1}(d_n^2)$ is the scalar $\|f_n\|_2^2$. The second term of the left-hand side of the above equivalence is called the *conditioned square function* of the martingale $(f_n)_{n\geq 1}$. Our main motivation comes from a remarkable article [16] by Junge and Xu where they extended the Burkholder/Rosenthal inequalities stated above in two directions. First, they found the right analogue of Burkholder inequalities for noncommutative martingales and second, the insight provided by the noncommutative case led to the right formulation of the corresponding inequality for the range $1 < p \leq 2$. To highlight the difference, we recall that in strong contrast with the classical case, (conditioned) square functions in the noncommutative case take several forms due to the row and column possibilities in the definition of martingale Hardy spaces and other related spaces. To motivate our consideration, recall that if $1 < p \leq 2$, and $x = (x_n)_{n\geq 1}$ is a noncommutative martingale, the norm on the (noncommutative) conditioned Hardy space $h^p$ introduced in [16] reads

$$\|x\|_{h^p} = \inf\left\{\left(\sum_{n\geq 1}\|a_n\|_p^p\right)^{1/p} + \left\|\left(\sum_{n\geq 1}\mathcal{E}_{n-1}(|b_n|^2)\right)^{1/2}\right\|_p + \left\|\left(\sum_{n\geq 1}\mathcal{E}_{n-1}(|c_n^*|^2)\right)^{1/2}\right\|_p\right\}$$

where the infimum runs over all decompositions $d_n = a_n + b_n + c_n$, with $a$, $b$ and $c$ being martingale difference sequences. With this norm, the noncommutative Burkholder/Rosenthal inequalities for $1 < p \leq 2$ from [16] can be formulated as follows:

(1.1) $$\|x\|_{h^p} \sim_{c_p} \sup_n \|x_n\|_p.$$

The crucial fact here is that for the case $1 < p \leq 2$, the $h_p$-norm requires that the given martingale be decomposed into three different martingales according to diagonal, column and row parts, respectively.

Inspired by (1.1) we consider the extremal case $p = 1$. Our main result (see Theorem 3.1) appears as decompositions of martingale difference sequences in the same spirit as in the definition of the $h^p$-norm when $1 < p < 2$. We can roughly state this as follows (see Theorem 3.1 for the full statement):

*there exists an absolute constant $K > 0$ such that if $(d_n)_{n\geq 1}$ is a martingale difference sequence in a noncommutative $L^2$-space associated with a finite von Neumann algebra $\mathcal{M}$, then there exists a decomposition $d_n = a_n + b_n + c_n$ satisfying the following weak-type $(1,1)$ inequality:*

$$(1.2) \quad \left\|\sum_{n\geq 1} a_n \otimes e_{n,n}\right\|_{L^{1,\infty}(\mathcal{M}\overline{\otimes}B(l^2))} + \left\|\left(\sum_{n\geq 1} \mathcal{E}_{n-1}(|b_n|^2)\right)^{1/2}\right\|_{1,\infty}$$
$$+ \left\|\left(\sum_{n\geq 1} \mathcal{E}_{n-1}(|c_n^*|^2)\right)^{1/2}\right\|_{1,\infty} \leq K \left\|\sum_{n\geq 1} d_n\right\|_1,$$

*where $(e_{i,j})_{i,j\geq 1}$ denotes the canonical matrix unit of $B(l^2)$.*

As in the weak-type $(1,1)$ inequality for square functions [31], our approach heavily depends on a noncommutative version of the classical Doob maximal inequality due to Cuculescu [3]. The proof is constructive and follows a line of ideas similar to that in [31]. The decomposition, however, has to be different since we have three separate terms as stated above.

Using general interpolation techniques and duality, our main result provides a new proof of the noncommutative analogue of Burkholder/Rosenthal inequalities from [16]. In fact, this approach improves considerably the estimates of the best constants from [16]. We obtain the optimal order of all the constants except for one case (see Theorem 4.1 below).

The paper is organized as follows: in Section 2, we briefly introduce the construction of noncommutative spaces and recall the general setup of martingales in noncommutative spaces along with formulations of various square functions that we will need throughout the paper. In Section 3, we formulate the appropriate weak-type $(1,1)$ inequality related to conditioned square functions. In Section 4, we revisit the noncommutative Burkholder/Rosenthal inequalities from [16]. Section 5 is devoted to study of sums of independent operators of mean zero in the sense of [16] and [35]. In the last section, we discuss some related results that point to some open problems.

**2. Preliminary definitions.** We use standard notation in operator algebras. We refer to [20, 33] for background on von Neumann algebra theory. Throughout all von Neumann algebras are assumed to be finite. Let $\mathcal{M}$ be a finite von Neumann algebra with a normal faithful finite trace $\tau$. The identity element of $\mathcal{M}$ is denoted by $\mathbf{1}$. For $1 \leq p \leq \infty$, we denote by $L^p(\mathcal{M}, \tau)$ [or simply $L^p(\mathcal{M})$] the noncommutative $L^p$-space associated with $(\mathcal{M}, \tau)$ (see, e.g., [6, 24]). Note that if $p = \infty$, we consider as customary $L^\infty(\mathcal{M}, \tau)$ as the von Neumann algebra $\mathcal{M}$ with the usual operator norm and recall that for $1 \leq p < \infty$, the norm on $L^p(\mathcal{M}, \tau)$ is defined by

$$\|x\|_p = (\tau(|x|^p))^{1/p}, \qquad x \in L^p(\mathcal{M}, \tau),$$

where $|x|$ is the usual modulus of $x$.

Assume that $\mathcal{M}$ is a $*$-subalgebra of $B(H)$ for a complex Hilbert space $H$. The elements of $L^p(\mathcal{M},\tau)$ can be viewed as closed densely defined operators on $H$. A closed densely defined operator $a$ on $H$ is said to be *affiliated with* $\mathcal{M}$ if $u^*au = a$ for all unitary $u$ in the commutant $\mathcal{M}'$ of $\mathcal{M}$. A closed densely defined operator $a$ on $H$ affiliated with $\mathcal{M}$ is said to be $\tau$-*measurable* if there exists $\lambda \geq 0$ such that $\tau(\chi_{(\lambda,\infty)}(|a|)) < \infty$ where $\chi_{(\lambda,\infty)}(|a|)$ denotes the spectral projection of $|a|$ corresponding to the characteristic function $\chi_{(\lambda,\infty)}(\cdot)$. For a measurable operator $a$, the generalized singular value function $\mu(a)$ is defined by

$$\mu_t(a) = \inf\{\lambda \geq 0 : \tau(\chi_{(\lambda,\infty)}(|a|)) \leq t\}, \qquad t \geq 0.$$

We refer to [9] for details and properties of the function $\mu(\cdot)$.

Of special interest in this paper is the noncommutative weak $L^1$-spaces associated with $(\mathcal{M},\tau)$ and denoted by $L^{1,\infty}(\mathcal{M},\tau)$. It is the collection of all $\tau$-measurable operators $x$ for which the quasi-norm

(2.1) $$\|x\|_{1,\infty} := \sup_{t>0} t\mu_t(x) = \sup_{\lambda>0} \lambda\tau(\chi_{(\lambda,\infty)}(|x|))$$

is finite. The following quasi-triangle inequality on elements of $L^{1,\infty}(\mathcal{M},\tau)$ holds and will be used repeatedly in the sequel. A short proof can be found in [31], Lemma 1.2.

LEMMA 2.1.  *For any $x_1, x_2$ in $L^{1,\infty}(\mathcal{M},\tau)$ and $\lambda > 0$,*

$$\lambda\tau(\chi_{(\lambda,\infty)}(|x_1+x_2|)) \leq 2\lambda\tau(\chi_{(\lambda/2,\infty)}(|x_1|)) + 2\lambda\tau(\chi_{(\lambda/2,\infty)}(|x_2|)).$$

For a complete, detailed and up-to-date presentation of noncommutative integration and noncommutative spaces, we refer to the recent survey [28].

Let us now recall the general setup for noncommutative martingales. Let $(\mathcal{M}_n)_{n\geq 1}$ be an increasing sequence of von Neumann subalgebras of $\mathcal{M}$ such that the union of $\mathcal{M}_n$'s is weak$^*$-dense in $\mathcal{M}$. For each $n \geq 1$, it is well known that there is a unique normal faithful conditional expectation $\mathcal{E}_n$ from $\mathcal{M}$ onto $\mathcal{M}_n$ such that $\tau \circ \mathcal{E}_n = \tau$. Moreover, $\mathcal{E}_n$ extends to a contractive projection from $L^p(\mathcal{M},\tau)$ onto $L^p(\mathcal{M}_n, \tau|_{\mathcal{M}_n})$ for every $1 \leq p < \infty$ which we will still denote by $\mathcal{E}_n$.

DEFINITION 2.2.  A noncommutative martingale with respect to the filtration $(\mathcal{M}_n)_{n\geq 1}$ is a sequence $x = (x_n)_{n\geq 1}$ in $L^1(\mathcal{M},\tau)$ such that

$$\mathcal{E}_n(x_{n+1}) = x_n \qquad \text{for all } n \geq 1.$$

If additionally, $x \in L^p(\mathcal{M}, \tau)$ for some $1 \leq p < \infty$, then $x$ is called an $L^p$-martingale. In this case, we set

$$\|x\|_p := \sup_{n \geq 1} \|x_n\|_p.$$

If $\|x\|_p < \infty$, then $x$ is called a bounded $L^p$-martingale. The difference sequence $dx = (dx_n)_{n \geq 1}$ of a martingale $x = (x_n)_{n \geq 1}$ is defined by

$$dx_n = x_n - x_{n-1}$$

with the usual convention that $x_0 = 0$. For concrete natural examples of noncommutative martingales, we refer to [27, 35].

We will now describe square functions of noncommutative martingales. Following [27], we will consider the following row and column versions of square functions: for a finite martingale $x = (x_n)_{n \geq 1}$, we denote by $dx$ the difference sequence as defined above. Set

$$S_C(x) = \left(\sum_{k \geq 1} |dx_k|^2\right)^{1/2} \quad \text{and} \quad S_R(x) = \left(\sum_{k \geq 1} |dx_k^*|^2\right)^{1/2}.$$

For $1 \leq p \leq \infty$ and any finite sequence $a = (a_n)_{n \geq 1}$ in $L^p(\mathcal{M}, \tau)$, set

$$\|a\|_{L^p(\mathcal{M};l_C^2)} = \left\|\left(\sum_{n \geq 1} |a_n|^2\right)^{1/2}\right\|_p, \qquad \|a\|_{L^p(\mathcal{M};l_R^2)} = \left\|\left(\sum_{n \geq 1} |a_n^*|^2\right)^{1/2}\right\|_p.$$

We recall the definitions of martingale Hardy spaces. Let $1 \leq p < \infty$; for a finite $L^p$-martingale $x$, set

$$\|x\|_{\mathcal{H}_C^p(\mathcal{M})} = \|dx\|_{L^p(\mathcal{M};l_C^2)} \quad \text{and} \quad \|x\|_{\mathcal{H}_R^p(\mathcal{M})} = \|dx\|_{L^p(\mathcal{M};l_R^2)}.$$

The space $\mathcal{H}_C^p(\mathcal{M})$ [resp. $\mathcal{H}_R^p(\mathcal{M})$] is defined as the completion of the collection of finite $L^p$-martingales under the norm $\|\cdot\|_{\mathcal{H}_C^p(\mathcal{M})}$ (resp. $\|\cdot\|_{\mathcal{H}_R^p(\mathcal{M})}$). The Hardy space of noncommutative martingales is defined as follows: if $1 \leq p < 2$,

$$\mathcal{H}^p(\mathcal{M}) = \mathcal{H}_C^p(\mathcal{M}) + \mathcal{H}_R^p(\mathcal{M})$$

equipped with the norm

$$\|x\|_{\mathcal{H}^p(\mathcal{M})} = \inf\{\|y\|_{\mathcal{H}_C^p(\mathcal{M})} + \|z\|_{\mathcal{H}_R^p(\mathcal{M})}\}$$

where the infimum runs over all pairs $(y, z) \in \mathcal{H}_C^p(\mathcal{M}) \times \mathcal{H}_R^p(\mathcal{M})$ such that $x = y + z$. For $2 \leq p < \infty$,

$$\mathcal{H}^p(\mathcal{M}) = \mathcal{H}_C^p(\mathcal{M}) \cap \mathcal{H}_R^p(\mathcal{M})$$

equipped with the intersection norm

$$\|x\|_{\mathcal{H}^p(\mathcal{M})} = \max\{\|x\|_{\mathcal{H}_C^p(\mathcal{M})}, \|x\|_{\mathcal{H}_R^p(\mathcal{M})}\}.$$

Below and throughout the rest of the paper we write $a_p \approx b_p$ as $p \to p_0$ to abbreviate the statement that there are two absolute positive constants $K_1$ and $K_2$ such that

$$K_1 \leq \frac{a_p}{b_p} \leq K_2 \qquad \text{for } p \text{ close to } p_0.$$

A fundamental result involving Hardy spaces is the noncommutative Burkholder–Gundy inequalities which we now state for further use.

THEOREM 2.3 (Noncommutative Burkholder–Gundy inequalities [19, 27, 31]). *Let $1 < p < \infty$ and let $x = (x_n)_{n=1}^\infty$ be an $L^p$-martingale. Then $x$ is bounded in $L^p(\mathcal{M}, \tau)$ if and only if $x$ belongs to $\mathcal{H}^p(\mathcal{M})$. If this is the case, then*

$$(\mathrm{BG}_p) \qquad \alpha_p^{-1} \|x\|_{\mathcal{H}^p(\mathcal{M})} \leq \|x\|_p \leq \beta_p \|x\|_{\mathcal{H}^p(\mathcal{M})}.$$

*Moreover, we have the following estimates for the best constants in $(\mathrm{BG}_p)$:*

   (i) $\alpha_p \approx (p-1)^{-1}$ as $p \to 1$;
   (ii) $\alpha_p \approx p$ as $p \to \infty$;
   (iii) $\beta_p \approx 1$ as $p \to 1$;
   (iv) $\beta_p \approx p$ as $p \to \infty$.

*These are the optimal orders of the constants $\alpha_p$ and $\beta_p$.*

The equivalence $(\mathrm{BG}_p)$ was first proved in the seminal paper [27]. The optimal orders of the constants involved follow from results in [19] and [31].

We now consider the conditioned versions of square functions and Hardy spaces developed in [15, 16].

Let $1 \leq p < \infty$. For a finite sequence $a = (a_n)_{n \geq 1}$ in $\mathcal{M}$, we define (recalling that $\mathcal{E}_0 = \mathcal{E}_1$)

$$\|a\|_{L^p_{\mathrm{cond}}(\mathcal{M}; l_C^2)} := \left\| \left( \sum_{n \geq 1} \mathcal{E}_{n-1}(a_n^* a_n) \right)^{1/2} \right\|_p.$$

It was shown in [15] that $\|\cdot\|_{L^p_{\mathrm{cond}}(\mathcal{M}; l_C^2)}$ is a norm on the vector space of all finite sequences in $\mathcal{M} \cap L^1(\mathcal{M}, \tau)$. The completion of the space of finite sequences in $\mathcal{M}$ equipped with the norm $\|\cdot\|_{L^p_{\mathrm{cond}}(\mathcal{M}; l_C^2)}$ will be denoted by $L^p_{\mathrm{cond}}(\mathcal{M}; l_C^2)$ and is the conditioned version of the space $L^p(\mathcal{M}; l_C^2)$ defined earlier. Similarly, we can define the conditioned row space $L^p_{\mathrm{cond}}(\mathcal{M}; l_R^2)$. A crucial fact that we will need in the sequel is that both spaces $L^p_{\mathrm{cond}}(\mathcal{M}; l_C^2)$ and $L^p_{\mathrm{cond}}(\mathcal{M}; l_R^2)$ can be realized as closed subspaces column and row (resp.) of the noncommutative space $L^p(\mathcal{M} \overline{\otimes} B(l^2(\mathbb{N}^2)))$ associated to the semifinite von Neumann algebra $\mathcal{M} \overline{\otimes} B(l^2(\mathbb{N}^2))$. For complete details on these facts we refer to [15].

Let $x = (x_n)_{n \geq 1}$ be a finite martingale in $L^2(\mathcal{M}, \tau)$; we set

$$\sigma_C(x) = \left( \sum_{n \geq 1} \mathcal{E}_{n-1}(|dx_n|^2) \right)^{1/2} \quad \text{and} \quad \sigma_R(x) = \left( \sum_{n \geq 1} \mathcal{E}_{n-1}(|dx_n^*|^2) \right)^{1/2}.$$

These will be called the column and row conditioned square functions, respectively. Observe that for $1 \leq p < \infty$,

$$\|\sigma_C(x)\|_p = \|dx\|_{L^p_{\text{cond}}(\mathcal{M}; l^2_C)} \quad \text{and} \quad \|\sigma_R(x)\|_p = \|dx\|_{L^p_{\text{cond}}(\mathcal{M}; l^2_R)}.$$

Let $h^p_C(\mathcal{M})$ [resp. $h^p_R(\mathcal{M})$] denote the closure in $L^p_{\text{cond}}(\mathcal{M}; l^2_C)$ [resp. $L^p_{\text{cond}}(\mathcal{M}; l^2_R)$] of all finite martingales in $\mathcal{M}$ (here we identified a martingale with its martingale difference sequence). Let $h^p_D(\mathcal{M})$ be the subspace of $l^p(L^p(\mathcal{M}, \tau))$ consisting of martingale difference sequences. Following [16], we define the conditioned version of martingale Hardy spaces as follows: if $1 \leq p < 2$,

$$h^p(\mathcal{M}) := h^p_D(\mathcal{M}) + h^p_C(\mathcal{M}) + h^p_R(\mathcal{M})$$

equipped with the norm

$$\|x\|_{h^p(\mathcal{M})} = \inf\{\|x^D\|_{h^p_D(\mathcal{M})} + \|x^C\|_{h^p_C(\mathcal{M})} + \|x^R\|_{h^p_R(\mathcal{M})}\},$$

where the infimum runs over all triples $(x^D, x^C, x^R) \in h^p_D(\mathcal{M}) \times h^p_C(\mathcal{M}) \times h^p_R(\mathcal{M})$ such that $x_n = x^D_n + x^C_n + x^R_n$ for all $n \geq 1$. For $2 \leq p < \infty$,

$$h^p(\mathcal{M}) := h^p_D(\mathcal{M}) \cap h^p_C(\mathcal{M}) \cap h^p_R(\mathcal{M})$$

equipped with the norm

$$\|x\|_{h^p(\mathcal{M})} = \max\{\|x\|_{h^p_D(\mathcal{M})}, \|x\|_{h^p_C(\mathcal{M})}, \|x\|_{h^p_R(\mathcal{M})}\}.$$

For $1 \leq p < \infty$, it is known from [16] that the linear space $h^p(\mathcal{M})$ is a Banach space.

**3. A weak-type inequality for conditioned square functions.** We will retain all notation introduced in the preliminaries. Unless specified otherwise, all adapted sequences are understood to be with respect to a fixed filtration of von Neumann subalgebras of $\mathcal{M}$. The principal result of this section is Theorem 3.1 below which can be viewed as a natural extension of the noncommutative Burkholder inequalities from [16] to the case $p = 1$.

THEOREM 3.1. *There is an absolute constant $K > 0$ such that if $x = (x_n)_{1 \leq n \leq N}$ is a finite $L^2$-bounded martingale, then there exist three adapted sequences $a = (a_n)_{1 \leq n \leq N}$, $b = (b_n)_{1 \leq n \leq N}$ and $c = (c_n)_{1 \leq n \leq N}$ in $L^2(\mathcal{M}, \tau)$ such that:*

($\alpha$) *for every $1 \leq n \leq N$, we have the decomposition*
$$dx_n = a_n + b_n + c_n;$$

($\beta$) *the $L^2$-norms satisfy*
$$\|a\|_{L^2(\mathcal{M},l_C^2)} + \|b\|_{L^2(\mathcal{M},l_C^2)} + \|c\|_{L^2(\mathcal{M},l_R^2)} \leq K\|x\|_2;$$

($\gamma$) *the conditioned square functions satisfy the weak-type inequality:*
$$\left\|\sum_{n=1}^{N} a_n \otimes e_{n,n}\right\|_{L^{1,\infty}(\mathcal{M}\overline{\otimes}B(l_N^2))} + \left\|\left(\sum_{n=1}^{N} \mathcal{E}_{n-1}(|b_n|^2)\right)^{1/2}\right\|_{1,\infty}$$
$$+ \left\|\left(\sum_{n=1}^{N} \mathcal{E}_{n-1}(|c_n^*|^2)\right)^{1/2}\right\|_{1,\infty} \leq K\|x\|_1,$$

where $(e_{i,j})_{1 \leq i,j \leq N}$ denotes the canonical matrix unit of $B(l_N^2)$.

As in previous weak-type results, our approach depends very heavily on a noncommutative version of the classical Doob weak-type $(1,1)$ maximal inequality, due to Cuculescu [3] (which we will recall below). We also note that through the standard decomposition of a general martingale into four positive martingales, the general case can be deduced easily from the special case of a positive martingale. Hence, without loss of generality, we can and do assume that the finite martingale $x = (x_n)_{1 \leq n \leq N}$ is a positive martingale and $\|x\|_1 = 1$.

3.1. *The decomposition of the martingale difference sequence.* We will explicitly describe the decomposition as stated. We start with the proposition (due to Cuculescu [3]) below which can be viewed as a substitute for the classical weak-type $(1,1)$ boundedness of maximal functions. We will state a version that incorporates the different properties that we need in the sequel. A short proof of the form stated below can be found in [29].

PROPOSITION 3.2 [3]. *For every $\lambda > 0$, there exists a finite sequence of decreasing projections $(q_n^{(\lambda)})_{1 \leq n \leq N}$ in $\mathcal{M}$ with:*

(a) *for every $1 \leq n \leq N$, $q_n^{(\lambda)} \in \mathcal{M}_n$;*
(b) $q_n^{(\lambda)} = q_{n-1}^{(\lambda)} \cdot \chi_{(0,\lambda]}(q_{n-1}^{(\lambda)} x_n q_{n-1}^{(\lambda)}) = \chi_{(0,\lambda]}(q_{n-1}^{(\lambda)} x_n q_{n-1}^{(\lambda)}) \cdot q_{n-1}^{(\lambda)}$. *In particular, $q_n^{(\lambda)}$ commutes with $q_{n-1}^{(\lambda)} x_n q_{n-1}^{(\lambda)}$;*
(c) $q_n^{(\lambda)} x_n q_n^{(\lambda)} \leq \lambda q_n^{(\lambda)}$;
(d) $(q_n^{(\lambda)})_{1 \leq n \leq N}$ *is a decreasing sequence and $\tau(\mathbf{1} - q_N^{(\lambda)}) \leq \lambda^{-1}$.*

The construction is based on the finite sequences $(q_n^{(2^k)})_{1\leq n \leq N}$ for $k \geq 1$. Following [30, 31], for $1 \leq n \leq N$, we set

$$p_{0,n} := \bigwedge_{k=0}^{\infty} q_n^{(2^k)}$$

(3.1)

$$p_{i,n} := \bigwedge_{k=i}^{\infty} q_n^{(2^k)} - \bigwedge_{k=i-1}^{\infty} q_n^{(2^k)} \qquad \text{for } i \geq 1.$$

Elementary but useful properties of the sequences $(p_{i,n})_{i\geq 0}$ that are relevant to our proof are collected in the following lemma.

LEMMA 3.3 ([30], Proposition 1.4). *For $1 \leq n \leq N$, the sequence of projections $(p_{i,n})_{i\geq 0}$ is pairwise disjoint with the following properties:*

(a) *for every $1 \leq n \leq N$ and $i \geq 0$, $p_{i,n} \in \mathcal{M}_n$;*
(b) *$\sum_{i=0}^{\infty} p_{i,n} = \mathbf{1}$ (for the strong operator topology);*
(c) *for every $m_0 \geq 0$, $\sum_{i=0}^{m_0} p_{i,n} \leq q_n^{(2^{m_0})}$.*

As in [30], we observe that for $n \geq 1$, $x = \sum_{j=0}^{\infty}\sum_{i=0}^{\infty} p_{i,n} x p_{j,n-1}$ for all $x \in L^1(\mathcal{M},\tau)$ (where the double sum may be taken using the $L^1$-norm). The (finite) sequences $a$, $b$ and $c$ are defined as follows:

(3.2) $\quad a_n := \begin{cases} 0, & \text{if } n=1, \\ \displaystyle\sum_{j=0}^{\infty}\sum_{i>j}(p_{i,n} - p_{i,n-1}p_{i,n})\,dx_n p_{j,n-1}, & \text{if } 2 \leq n \leq N; \end{cases}$

(3.3) $\quad b_n := \begin{cases} \displaystyle\sum_{j=0}^{\infty}\sum_{i\leq j} p_{i,1}\,dx_1 p_{j,1}, & \text{if } n=1, \\ \displaystyle\sum_{j=0}^{\infty}\sum_{i\leq j} p_{i,n}\,dx_n p_{j,n-1}, & \text{if } 2 \leq n \leq N; \end{cases}$

(3.4) $\quad c_n := \begin{cases} \displaystyle\sum_{j=0}^{\infty}\sum_{i>j} p_{i,1}\,dx_1 p_{j,1}, & \text{if } n=1, \\ \displaystyle\sum_{j=0}^{\infty}\sum_{i>j} p_{i,n-1}p_{i,n}\,dx_n p_{j,n-1}, & \text{if } 2 \leq n \leq N. \end{cases}$

It is clear from this construction that for every $1 \leq n \leq N$, $a_n$, $b_n$ and $c_n$ belong to $L^2(\mathcal{M}_n, \tau|_{\mathcal{M}_n})$ and $dx_n = a_n + b_n + c_n$. Moreover, using boundedness of the triangular truncations in $L^2(\mathcal{M},\tau)$, it is straightforward that condition $(\beta)$ of the theorem is satisfied. Thus, it remains to prove the weak-type $(1,1)$ inequality as stated in condition $(\gamma)$ of Theorem 3.1.

3.2. *Proof of the weak-type $(1,1)$ inequalities.* The proof is separated into two parts highlighted in Propositions A and B below.

PROPOSITION A. *There exists an absolute constant $\kappa > 0$ such that*

$$\left\| \sum_{n \geq 1} a_n \otimes e_{n,n} \right\|_{L^{1,\infty}(\mathcal{M} \overline{\otimes} B(l_N^2))} \leq \kappa.$$

We will start by recording some basic lemmas for further use in the proof. For a given operator $x \in \mathcal{M}$, we denote by $l(x)$ [resp. $r(x)$] the left (resp. right) support projection of $x$ (see, e.g., [33], page 134, for definitions). We need the following observation:

LEMMA 3.4. *Let $n \geq 2$ and $i \geq 1$. Then:*

(i) $r(p_{i,n} - p_{i,n-1}p_{i,n}) \leq p_{i,n}$;
(ii) $l(p_{i,n} - p_{i,n-1}p_{i,n}) \leq \bigwedge_{k=i-1}^{\infty} q_{n-1}^{(2^k)} - \bigwedge_{k=i-1}^{\infty} q_n^{(2^k)}$.

PROOF. The first statement is trivial from the definition of right support projections. For the left support projections, we observe from the definition of $p_{i,n}$'s and the fact that $(q_n^{(2^k)})_{n \geq 1}$ is decreasing that for $n \geq 2$ and $i \geq 1$,

$$p_{i,n} - p_{i,n-1}p_{i,n} = \bigwedge_{k=i-1}^{\infty} q_{n-1}^{(2^k)} \bigwedge_{k=i}^{\infty} q_n^{(2^k)} - \bigwedge_{k=i-1}^{\infty} q_n^{(2^k)}.$$

The statement then follows directly from the definition of left support projections. □

LEMMA 3.5. *For $n \geq 2$ and $i \geq 1$, set $r_{i,n} = r(p_{i,n} - p_{i,n-1}p_{i,n})$; then for every $m_0 \in \mathbb{N}$,*

$$\sum_{n \geq 2} \tau \left( \sum_{i \geq m_0+1} r_{i,n} \right) \leq 4.2^{-m_0}.$$

PROOF. First, note that for every $n \geq 2$ and $i \geq 1$, $r_{i,n}$ is equivalent to $l(p_{i,n} - p_{i,n-1}p_{i,n})$ (see, e.g., [33], Proposition 1.5, page 292) and therefore from Lemma 3.4, we get that $\tau(r_{n,i}) = \tau(l(p_{i,n}-p_{i,n-1}p_{i,n})) \leq \tau(\bigwedge_{k=i-1}^{\infty} q_{n-1}^{(2^k)} - \bigwedge_{k=i-1}^{\infty} q_n^{(2^k)})$. Hence, from this estimate, we deduce that

$$\sum_{n \geq 2} \tau \left( \sum_{i \geq m_0+1} r_{i,n} \right) \leq \sum_{n \geq 2} \tau \left( \sum_{i \geq m_0+1} \bigwedge_{k=i-1}^{\infty} q_{n-1}^{(2^k)} - \bigwedge_{k=i-1}^{\infty} q_n^{(2^k)} \right)$$

$$= \sum_{i \geq m_0+1} \tau \left( \sum_{n \geq 2} \bigwedge_{k=i-1}^{\infty} q_{n-1}^{(2^k)} - \bigwedge_{k=i-1}^{\infty} q_n^{(2^k)} \right)$$

$$= \sum_{i \geq m_0+1} \tau\left(\bigwedge_{k=i-1}^{\infty} q_1^{(2^k)} - \bigwedge_{k=i-1}^{\infty} q^{(2^k)}\right)$$

$$\leq \sum_{i \geq m_0+1} \tau\left(\mathbf{1} - \bigwedge_{k=i-1}^{\infty} q^{(2^k)}\right)$$

$$= \sum_{i \geq m_0+1} \tau\left(\bigvee_{k=i-1}^{\infty} (\mathbf{1} - q^{(2^k)})\right)$$

$$\leq \sum_{i \geq m_0+1} \sum_{k \geq i-1} \tau(\mathbf{1} - q^{(2^k)}).$$

Thus, we conclude from Proposition 3.2 that

$$\sum_{n \geq 2} \tau\left(\sum_{i \geq m_0+1} r_{i,n}\right) \leq \sum_{i \geq m_0+1}\left(\sum_{k \geq i-1} 2^{-k}\right) = 4 \cdot 2^{-m_0},$$

which proves the lemma. □

LEMMA 3.6. *For $n \geq 2$, let $h_n = \sum_{i \geq 0} p_{i,n} - p_{i,n-1}p_{i,n} \in \mathcal{M}_n$ and set*

$$h := \sum_{n \geq 2} h_n \otimes e_{n,n} \in \mathcal{M} \overline{\otimes} B(l_N^2).$$

*Then $\max\{\|h\|_\infty, \|h\|_2\} \leq 2$.*

PROOF. Note that $\|h\|_\infty = \sup_{n \geq 2} \|h_n\|_\infty$ and for every $n \geq 2$, $h_n = \mathbf{1} - \sum_{i \geq 1} p_{i,n-1}p_{i,n}$. Since $|\sum_{i \geq 1} p_{i,n-1}p_{i,n}|^2 = \sum_{i \geq 1} p_{i,n}p_{i,n-1}p_{i,n} \leq \mathbf{1}$, the first assertion follows. For the $L^2$-norm, it is clear that $\|h_n\|_2^2 = \sum_{i \geq 1} \|p_{i,n} - p_{i,n-1}p_{i,n}\|_2^2 \leq \sum_{i \geq 1} \tau(r_{i,n})$. From Lemma 3.5, we deduce that

$$\|h\|_2^2 \leq \sum_{n \geq 2} \sum_{i \geq 1} \tau(r_{i,n}) \leq 4$$

which shows the desired estimate on the $L^2$-norm. □

We are now ready to provide the proof of Proposition A.

PROOF OF PROPOSITION A. We denote by tr the usual trace of $B(l_N^2)$. From the definition of $\|\cdot\|_{1,\infty}$, it is enough to show the existence of a numerical constant $\kappa > 0$ such that for every $\lambda > 0$,

(3.5) $$\tau \otimes \mathrm{tr}\left(\chi_{(\lambda,\infty)}\left(\left|\sum_{n \geq 2} a_n \otimes e_{n,n}\right|\right)\right) \leq \kappa \lambda^{-1}.$$

Since the trace $\tau \otimes \mathrm{tr}$ is not normalized, we have to verify (3.5) for the full range $0 < \lambda < \infty$. We separate the proof into two separate cases according to $\lambda \geq 1$ or $0 < \lambda < 1$.

*Case* 1. Assume that $\lambda \geq 1$. For this case, it is enough to verify (3.5) for $\lambda = 2^{m_0}$ for $m_0 \geq 0$. To simplify the notation, we set $\Theta := \sum_{n \geq 2} a_n \otimes e_{n,n}$.

We first observe that

$$\Theta = h \cdot \sum_{n \geq 2} \left( \sum_{j \geq 0} \sum_{i > j} r_{i,n} dx_n p_{j,n-1} \right) \otimes e_{n,n}$$
$$= h \cdot \sum_{n \geq 2} \left( \sum_{j \geq 0} \sum_{i > j} (r_{i,n} \otimes e_{n,n}) \cdot (dx_n \otimes e_{n,n}) \cdot (p_{j,n-1} \otimes e_{n,n}) \right).$$

Consider the following projection in $\mathcal{M} \overline{\otimes} B(l_N^2)$:

(3.6) $$v_0 := \sum_{n \geq 2} \left( \sum_{i \geq m_0 + 1} r_{i,n} \right) \otimes e_{n,n}.$$

Then from the estimate in Lemma 3.5, we have $\tau \otimes \operatorname{tr}(v_0) \leq 4\lambda^{-1}$. Moreover, we can write $\Theta$ using the projection $v_0$ as

$$\Theta = h \cdot \sum_{n \geq 2} \left( \sum_{j=0}^{m_0} \sum_{j < i \leq m_0} r_{i,n} dx_n p_{j,n-1} \right) \otimes e_{n,n} + h \cdot v_0 \cdot \Pi,$$

where $\Pi = \sum_{n \geq 2} (\sum_{j \geq 0} \sum_{i > \max(j,m_0)} r_{i,n} dx_n p_{j,n-1}) \otimes e_{n,n}$. We can split the trace according to Lemma 2.1 and get

$$\tau \otimes \operatorname{tr}(\chi_{(\lambda,\infty)}(|\Theta|)) \leq 2\tau \otimes \operatorname{tr}(\chi_{(\lambda/2,\infty)}(|h \cdot \Gamma|)) + 2\tau \otimes \operatorname{tr}(\chi_{(\lambda/2,\infty)}(|h \cdot v_0 \cdot \Pi|))$$

where $\Gamma = \sum_{n \geq 2} (\sum_{j=0}^{m_0} \sum_{j < i \leq m_0} r_{i,n} dx_n p_{j,n-1}) \otimes e_{n,n}$. A fortiori,

(3.7) $$\tau \otimes \operatorname{tr}(\chi_{(\lambda,\infty)}(|\Theta|)) \leq 8\|h\|_\infty^2 \lambda^{-2} \|\Gamma\|_2^2$$
$$+ 2\tau \otimes \operatorname{tr}(\chi_{(\lambda^2/4,\infty)}(hv_0 \Pi \Pi^* v_0 h^*)).$$

Observe (see, e.g., [9]) that

$$\tau \otimes \operatorname{tr}(\chi_{(\lambda^2/4,\infty)}(hv_0 \Pi \Pi^* v_0 h^*))$$
$$= \int_0^\infty \chi_{(\lambda^2/4,\infty)}\{\mu_t(hv_0 \Pi \Pi^* v_0 h^*)\} dt$$
$$\leq \int_0^\infty \chi_{(\lambda^2/4\|h\|_\infty^2,\infty)}\{\mu_t(v_0 \Pi \Pi^* v_0)\} dt$$

where the singular value $\mu_t(\cdot)$ is relative to $\mathcal{M} \overline{\otimes} B(l_N^2)$. Therefore,

$$\tau \otimes \operatorname{tr}(\chi_{(\lambda^2/4,\infty)}(hv_0 \Pi \Pi^* v_0 h^*)) \leq \tau \otimes \operatorname{tr}(v_0) \leq 4\lambda^{-1}.$$

Moreover, it is clear that for every $s \geq 1$, $\sum_{i=0}^{m_0} r_{i,s} \leq \sum_{i=0}^{m_0} p_{i,s} \leq q_s^{(\lambda)}$ and thus $\Gamma = \sum_{n \geq 2} (\sum_{j=0}^{m_0} \sum_{j < i \leq m_0} r_{i,n} (q_n^{(\lambda)} dx_n q_{n-1}^{(\lambda)}) p_{j,n-1}) \otimes e_{n,n}$. From the $L^2$-boundedness of triangular truncations, it follows that

$$\|\Gamma\|_2^2 \leq \sum_{n \geq 2} \|q_n^{(\lambda)} dx_n q_{n-1}^{(\lambda)}\|_2^2.$$

Combining the above estimates, (3.7) implies

$$(3.8) \quad \tau \otimes \mathrm{tr}(\chi_{(\lambda,\infty)}(|\Theta|)) \leq 32\lambda^{-2} \sum_{n\geq 2} \|q_n^{(\lambda)} dx_n q_{n-1}^{(\lambda)}\|_2^2 + 8\lambda^{-1}.$$

Therefore it remains to estimate $\|q_n^{(\lambda)} dx_n q_{n-1}^{(\lambda)}\|_2^2$ for $n \geq 2$. This is done in the next lemma.

LEMMA 3.7. *For every* $2 \leq n \leq N$, $\|q_n^{(\lambda)} dx_n q_{n-1}^{(\lambda)}\|_2 \leq \|q_n^{(\lambda)} x_n q_n^{(\lambda)} - q_{n-1}^{(\lambda)} x_{n-1} q_{n-1}^{(\lambda)}\|_2$ *and*

$$\|q_1^{(\lambda)} x_1 q_1^{(\lambda)}\|_2^2 + \sum_{n=2}^{N} \|q_n^{(\lambda)} x_n q_n^{(\lambda)} - q_{n-1}^{(\lambda)} x_{n-1} q_{n-1}^{(\lambda)}\|_2^2 \leq 2\lambda.$$

PROOF. We will simply write $q_n$ for $q_n^{(\lambda)}$. For every $2 \leq n \leq N$, we note that $q_n \leq q_{n-1}$ and $q_n$ commutes with $q_{n-1} x_n q_{n-1}$ [Proposition 3.2(c)]. Therefore,

$$\begin{aligned}
\|q_n dx_n q_{n-1}\|_2^2 &= \tau(q_n\, dx_n q_{n-1}\, dx_n q_n) \\
&= \tau(q_n(x_n - x_{n-1}) q_{n-1} (x_n - x_{n-1}) q_n) \\
&= \tau(|[q_{n-1} x_n q_{n-1} - q_{n-1} x_{n-1} q_{n-1}] q_n|^2) \\
&= \tau(|[q_n x_n q_n - q_{n-1} x_{n-1} q_{n-1}] q_n|^2) \\
&\leq \|q_n x_n q_n - q_{n-1} x_{n-1} q_{n-1}\|_2^2.
\end{aligned}$$

The proof of the second inequality can be found in [31], Proposition E. □

We can now conclude by combining (3.8) and Lemma 3.7 that

$$\tau \otimes \mathrm{tr}(\chi_{(\lambda,\infty)}(|\Theta|)) \leq 72\lambda^{-1}.$$

Thus the proof for the case $\lambda \geq 1$ is complete.

*Case* 2. Assume $0 < \lambda < 1$. Let $m_0 = 1$ and let $v_0$ be the corresponding projection in $\mathcal{M} \overline{\otimes} B(l_N^2)$ as in (3.5) above. Then (as in the previous case), we write

$$\Theta = h \cdot v_0 \cdot \Pi + h \cdot \left( \sum_{n \geq 2} r_{1,n} q_n^{(2)} dx_n q_{n-1}^{(2)} p_{0 \cdot n-1} \otimes e_{n,n} \right).$$

As $\tau \otimes \mathrm{tr}(v_0) \leq 2$, the argument in the previous case gives

$\tau \otimes \mathrm{tr}(\chi_{(\lambda,\infty)}(|\Theta|))$

$\leq 2\tau \otimes \mathrm{tr}\left( \chi_{(\lambda/2,\infty)}\left( \left| h \cdot \left( \sum_{n\geq 2} r_{1,n} q_n^{(2)} dx_n q_{n-1}^{(2)} p_{0\cdot n-1} \otimes e_{n,n} \right) \right| \right) \right) + 4$

$$\leq 4\lambda^{-1} \left\| h \cdot \left( \sum_{n \geq 2} r_{1,n} q_n^{(2)} \, dx_n q_{n-1}^{(2)} p_{0 \cdot n-1} \otimes e_{n,n} \right) \right\|_1 + 4$$

$$= 4\lambda^{-1} \sum_{n \geq 2} \| h_n r_{1,n} q_n^{(2)} \, dx_n q_{n-1}^{(2)} p_{0 \cdot n-1} \|_1 + 4$$

$$\leq 4\lambda^{-1} \sum_{n \geq 2} \| h_n \|_2 \| q_n^{(2)} \, dx_n q_{n-1}^{(2)} \|_2 + 4.$$

Since $\|h\|_2 = (\sum_{n \geq 2} \|h_n\|_2^2)^{1/2} \leq 2$, Hölder's inequality implies

$$\tau \otimes \operatorname{tr}(\chi_{(\lambda,\infty)}(|\Theta|)) \leq 8\lambda^{-1} \left( \sum_{n \geq 2} \| q_n^{(2)} \, dx_n q_{n-1}^{(2)} \|_2^2 \right)^{1/2} + 4.$$

We now apply Lemma 3.7 to conclude that

$$\tau \otimes \operatorname{tr}(\chi_{(\lambda,\infty)}(|\Theta|)) \leq 32\lambda^{-1} + 4 \leq 36\lambda^{-1}.$$

This proves (3.5) for the case $0 < \lambda < 1$ and combined with the previous case, the proof of Proposition A is complete. □

PROPOSITION B. *There exists an absolute constant $\kappa > 0$ so that*

$$\left\| \left( \sum_{n \geq 1} \mathcal{E}_{n-1}(|b_n|^2) \right)^{1/2} \right\|_{1,\infty} + \left\| \left( \sum_{n \geq 1} \mathcal{E}_{n-1}(|c_n^*|^2) \right)^{1/2} \right\|_{1,\infty} \leq \kappa.$$

PROOF. We begin by highlighting the forms of the conditioned square functions relative to the sequences $b$ and $c$. The proof of the following lemma is just a notational adjustment of [30], Lemma 2.2 and is left to the interested reader.

LEMMA 3.8. *For the sequences defined above, we have:*

(a) $|b_1|^2 = \sum_{l=0}^{\infty} \sum_{j=0}^{\infty} \sum_{i \leq \min(l,j)} p_{l,1} \, dx_1 p_{i,1} \, dx_1 p_{j,1}$;
(b) $\mathcal{E}_{n-1}(|b_n|^2) = \sum_{l=0}^{\infty} \sum_{j=0}^{\infty} \sum_{i \leq \min(l,j)} p_{l,n-1} \mathcal{E}_{n-1}[dx_n p_{i,n} \, dx_n] p_{j,n-1}$ *for* $2 \leq n \leq N$;
(c) $|c_1^*|^2 = \sum_{l=1}^{\infty} \sum_{j=1}^{\infty} \sum_{i < \min(l,j)} p_{l,1} \, dx_1 p_{i,1} \, dx_1 p_{j,1}$;
(d) *for* $2 \leq n \leq N$,

$$\mathcal{E}_{n-1}(|c_n^*|^2)$$
$$= \sum_{l=1}^{\infty} \sum_{j=1}^{\infty} \sum_{i < \min(l,j)} p_{l,n-1} \mathcal{E}_{n-1}[p_{l,n} dx_n p_{i,n-1} \, dx_n p_{j,n}] p_{j,n-1},$$

*where the sums are taken using the $L^1$-norm.*

From the preceding lemma, we remark that $\mathcal{E}_{n-1}(|b_n|^2)$ and $\mathcal{E}_{n-1}(|c_n^*|^2)$ are essentially of the same form so we will provide the complete detail on the appropriate estimate of $\|(\sum_{n\geq 1} \mathcal{E}_{n-1}(|b_n|^2))^{1/2}\|_{1,\infty}$ and only point out the (minimal) adjustment needed for $\|(\sum_{n\geq 1} \mathcal{E}_{n-1}(|c_n^*|^2))^{1/2}\|_{1,\infty}$. In particular, we need to verify the existence of an absolute constant $\kappa > 0$ such that

$$\left\| \left( \sum_{n\geq 1} \mathcal{E}_{n-1}(|b_n|^2) \right)^{1/2} \right\|_{1,\infty} \leq \kappa. \tag{3.9}$$

As above, this is equivalent to showing the existence of a numerical constant $\kappa > 0$ such that for every $\lambda > 0$,

$$\tau\left( \chi_{(\lambda,\infty)}\left( \left( \sum_{n\geq 1} \mathcal{E}_{n-1}(|b_n|^2) \right)^{1/2} \right) \right) \leq \kappa \lambda^{-1}. \tag{3.10}$$

We note that since $\tau$ is normalized, it is enough to verify the existence of such constant for $\lambda = 2^{m_0}$ where $m_0 \geq 0$. The proof basically follows the steps used in the previous proposition. Throughout the proof, let $\sigma_C := (\sum_{n\geq 1} \mathcal{E}_{n-1}(|b_n|^2))^{1/2}$.

Consider the projection

$$w_0 := \sum_{i=0}^{m_0} p_{i,N} = \bigwedge_{k=m_0}^{\infty} q_N^{(2^k)} \tag{3.11}$$

and

$$\gamma_0 = \left| \sum_{j=0}^{m_0} \sum_{i\leq j} p_{i,1}\, dx_1 p_{j,1} \right|^2 \tag{3.12}$$
$$+ \sum_{n=2}^{N} \sum_{l=0}^{m_0} \sum_{j=0}^{m_0} \sum_{i\leq \min(l,j)} p_{l,n-1} \mathcal{E}_{n-1}[dx_n p_{i,n}\, dx_n] p_{j,n-1}.$$

It is clear that $\tau(\mathbf{1} - w_0) \leq \sum_{k\geq m_0+1} \tau(\mathbf{1} - q_N^{(2^k)}) \leq 2^{-m_0} = \lambda^{-1}$. Moreover, since for every $1 \leq s \leq N$, $w_0 \leq \sum_{i=1}^{m_0} p_{i,s}$, we have $w_0 \sigma_C^2 w_0 = w_0 \gamma_0 w_0$. We deduce from Lemma 2.1 that

$$\tau(\chi_{(\lambda,\infty)}(\sigma_C)) \leq 2\tau(\chi_{(\lambda/2,\infty)}(|\sigma_C w_0|)) + 2\tau(\chi_{(\lambda/2,\infty)}(|\sigma_C(\mathbf{1}-w_0)|))$$
$$\leq 2\tau(\chi_{(\lambda^2/4,\infty)}(w_0 \gamma_0 w_0)) + 2\lambda^{-1}$$
$$\leq 8\lambda^{-2}\|w_0 \gamma_0 w_0\|_1 + 2\lambda^{-1}.$$

We remark that since the expectations are $\tau$-invariant,

$$\|w_0 \gamma_0 w_0\|_1 \leq \left\| \sum_{j=0}^{m_0} \sum_{i\leq j} p_{i,1}\, dx_1 p_{j,1} \right\|_2^2$$

$$+ \tau \bigg( \sum_{n=2}^{N} \sum_{l=0}^{m_0} \sum_{j=0}^{m_0} \sum_{i \leq \min(l,j)} p_{l,n-1} \mathcal{E}_{n-1}[dx_n p_{i,n} \, dx_n] p_{j,n-1} \bigg)$$

$$= \bigg\| \sum_{j=0}^{m_0} \sum_{i \leq j} p_{i,1} \, dx_1 p_{j,1} \bigg\|_2^2$$

$$+ \tau \bigg( \sum_{n=2}^{N} \sum_{j=0}^{m_0} \sum_{i \leq j} p_{j,n-1} \, dx_n p_{i,n} \, dx_n p_{j,n-1} \bigg).$$

Moreover, from the fact that $\sum_{i=1}^{m_0} p_{i,s} \leq q_s^{(\lambda)}$ when $1 \leq s \leq m_0$, we have

$$\|w_0 \gamma_0 w_0\|_1 \leq \bigg\| \sum_{j=0}^{m_0} \sum_{i \leq j} p_{i,1}(q_1^{(\lambda)} \, dx_1 q_1^{(\lambda)}) p_{j,1} \bigg\|_2^2$$

$$+ \tau \bigg( \sum_{n=2}^{N} \sum_{j=0}^{m_0} p_{j,n-1}(q_{n-1}^{(\lambda)} \, dx_n q_n^{(\lambda)} \, dx_n q_{n-1}^{(\lambda)}) p_{j,n-1} \bigg)$$

$$\leq \|q_1^{(\lambda)} \, dx_1 q_1^{(\lambda)}\|_2^2 + \sum_{n=2}^{N} \|q_n^{(\lambda)} \, dx_n q_{n-1}^{(\lambda)}\|_2^2.$$

We can now apply Lemma 3.7 to conclude that

(3.13) $$\tau(\chi_{(\lambda,\infty)}(\sigma_C)) \leq 18 \lambda^{-1}.$$

For the second part of the proposition, let $\sigma_R := (\sum_{n \geq 1} \mathcal{E}_{n-1}(|c_n^*|^2))^{1/2}$. A notational adjustment of the argument used above leads to the inequality

$$\tau(\chi_{(\lambda,\infty)}(\sigma_R)) \leq 8 \lambda^{-2} \Xi + 2 \lambda^{-1}$$

with

$$\Xi \leq \|q_1^{(\lambda)} \, dx_1 q_1^{(\lambda)}\|_2^2 + \tau \bigg( \sum_{n=2}^{N} \sum_{j=0}^{m_0} \sum_{i \leq j} p_{j,n-1}[p_{j,n} \, dx_n p_{i,n-1} \, dx_n p_{j,n}] p_{j,n-1} \bigg)$$

$$\leq \|q_1^{(\lambda)} \, dx_1 q_1^{(\lambda)}\|_2^2 + \tau \bigg( \sum_{n=2}^{N} \sum_{j=0}^{m_0} p_{j,n-1}[p_{j,n} q_n^{(\lambda)} \, dx_n q_{n-1}^{(\lambda)} \, dx_n q_n^{(\lambda)} p_{j,n}] p_{j,n-1} \bigg)$$

$$\leq \|q_1^{(\lambda)} \, dx_1 q_1^{(\lambda)}\|_2^2 + \sum_{n=2}^{N} \|q_{n-1}^{(\lambda)} \, dx_n q_n^{(\lambda)}\|_2^2,$$

and thus as above, $\Xi \leq 2\lambda$ and therefore it follows that

(3.14) $$\tau(\chi_{(\lambda,\infty)}(\sigma_R)) \leq 18 \lambda^{-1}.$$

Proposition B follows from combining (3.13) and (3.14). □

The weak-type $(1,1)$ inequality in Theorem 3.1 clearly follows by combining Propositions A and B. The proof is complete.

REMARK 3.9. (i) In the proofs of Propositions A and B above, one can also use the noncommutative analogue of Gundy's decomposition recently obtained in [26]. This, however, does not lead to substantial simplification of the more traditional approach used above.

(ii) We note that the decomposition in Theorem 3.1 is only with adapted sequences. We do not know if such decomposition can be accomplished with martingale difference sequences. From the Cauchy–Schwarz inequality on conditional expectations, it follows that Proposition B is still valid with $(b_n - \mathcal{E}_{n-1}(b_n))_{n \geq 1}$ and $(c_n - \mathcal{E}_{n-1}(c_n))_{n \geq 1}$ in place of $(b_n)_{n \geq 1}$ and $(c_n)_{n \geq 1}$, respectively. At this time we do not know if Proposition A holds with $(a_n - \mathcal{E}_{n-1}(a_n))_{n \geq 1}$. This, however, will not affect the argument used in the main application in the next section. See Proposition 3.10 below for a particular case where the decomposition can be done with martingale difference sequences.

After the first draft of this paper was written, we learned that J. Parcet had proved in [25] two nonequivalent weak-type $(1,1)$ inequalities analogous to Burkholder inequalities (see [25], Theorem A, Theorem B, Corollary C). His paper, however, considered only classical martingales on probability spaces. It turns out that Theorem 3.1 above when applied to classical martingales gives the same result as [25], Corollary C. His primary tool is the classical Davis decomposition. We also obtain a noncommutative analogue of [25], Theorem B under similar assumption. Following [25] for classical martingales, we say that a positive noncommutative martingale $x = (x_n)_{n \geq 1}$ is $k$-regular (for some constant $k > 1$) if for every $n \geq 2$,

$$x_n \leq k x_{n-1}.$$

PROPOSITION 3.10. *There is an absolute constant $C > 0$ such that if $x = (x_n)_{n \geq 1}$ is a $k$-regular $L^2$-bounded martingale, then we can decompose $x = y + z$ as sum of two martingales satisfying the inequality*

$$\|\sigma_C(y)\|_{1,\infty} + \|\sigma_R(z)\|_{1,\infty} \leq C k^2 \|x\|_1.$$

It should be noted that in strong contrast with the general case, the conditioned Hardy norms for $k$-regular martingales do not require the diagonal term. This difference was already observed by Parcet in [25] for classical martingales. We also refer the reader to [10], pages 124–127, for the case of (commutative) predictable martingale $f$ with $\sup|f_n| \in L^1$. Our proof below has to be different from the classical case since as expected we have to take into account the decomposition into row and column parts. The decomposition, however, is exactly the same as the one used in the case of square functions from [31].

SKETCH OF THE PROOF OF PROPOSITION 3.10.    Define the martingales $y = (y_n)_{n\geq 1}$ and $z = (z_n)_{n\geq 1}$ exactly as in [31]:

(3.15)
$$dy_1 := \sum_{j=0}^{\infty}\sum_{i\leq j} p_{i,1}\, dx_1 p_{j,1},$$
$$dy_n := \sum_{j=0}^{\infty}\sum_{i\leq j} p_{i,n-1}\, dx_n p_{j,n-1} \qquad \text{for } n \geq 2;$$

and

(3.16)
$$dz_1 := \sum_{j=0}^{\infty}\sum_{i>j} p_{i,1}\, dx_1 p_{j,1},$$
$$dz_n := \sum_{j=0}^{\infty}\sum_{i>j} p_{i,n-1}\, dx_n p_{j,n-1} \qquad \text{for } n \geq 2.$$

We refer to [31] for the fact that $(dy_n)_{n\geq 1}$ and $(dz_n)_{n\geq 1}$ are martingale difference sequences and satisfy $dx_n = dy_n + dz_n$ for $n \geq 1$. The proof on the estimate of $\|\sigma_C(y)\|_{1,\infty} + \|\sigma_R(z)\|_{1,\infty}$ is quite elementary and essentially follows the steps used in Proposition B above so we will only sketch the main points.

As in the proof of Proposition B, it suffices to estimate $\tau(\chi_{(\lambda,\infty)}(\sigma_C(y)))$ for dyadic $\lambda = 2^{m_0}$ where $m_0 \in \mathbb{N}$. We remark first that for $n \geq 2$,

$$\mathcal{E}_{n-1}(|b_n|^2) = \sum_{l=0}^{\infty}\sum_{j=0}^{\infty}\sum_{i\leq \min(l,j)} p_{l,n-1}\mathcal{E}_{n-1}[dx_n p_{i,n-1}\, dx_n]p_{j,n-1}.$$

If $w_0$ is the projection defined in (3.11), then

$$\tau(\chi_{(\lambda,\infty)}(\sigma_C(y))) \leq 2\tau(\chi_{(\lambda/2,\infty)}(|\sigma_C(y)w_0|)) + 2\tau(\mathbf{1} - w_0)$$
$$\leq 8\lambda^{-2}\|\sigma_C(y)w_0\|_2^2 + 2\lambda^{-1}.$$

A just notational adjustment of the argument used in the proof of Proposition B leads to

$$\|\sigma_C(y)w_0\|_2^2 \leq \|q_1^{(\lambda)}\, dx_1 q_1^{(\lambda)}\|_2^2 + \sum_{n\geq 2} \|q_{n-1}^{(\lambda)}\, dx_n q_{n-1}^{(\lambda)}\|_2^2.$$

Splitting the quantity $q_{n-1}^{(\lambda)} dx_n q_{n-1}^{(\lambda)}$ to $q_{n-1}^{(\lambda)} dx_n q_n^{(\lambda)} + q_{n-1}^{(\lambda)} dx_n (q_{n-1}^{(\lambda)} - q_n^{(\lambda)})$ for all $n \geq 2$ gives

$$\|\sigma_C(y) w_0\|_2^2 \leq \left( \|q_1^{(\lambda)} dx_1 q_1^{(\lambda)}\|_2^2 + 2 \sum_{n \geq 2} \|q_{n-1}^{(\lambda)} dx_n q_n^{(\lambda)}\|_2^2 \right)$$

$$+ 2 \sum_{n \geq 2} \|q_{n-1}^{(\lambda)} dx_n q_{n-1}^{(\lambda)} (q_{n-1}^{(\lambda)} - q_n^{(\lambda)})\|_2^2$$

$$= I + II.$$

It follows from Lemma 3.7 that $I \leq 4\lambda$. For $II$, it is immediate that

$$II \leq 2 \sum_{n \geq 2} \|q_{n-1}^{(\lambda)} dx_n q_{n-1}^{(\lambda)}\|_\infty^2 \tau(q_{n-1}^{(\lambda)} - q_n^{(\lambda)}).$$

Since $x$ is $k$-regular, we deduce for all $n \geq 2$ that

$$\|q_{n-1}^{(\lambda)} dx_n q_{n-1}^{(\lambda)}\|_\infty \leq \|q_{n-1}^{(\lambda)} x_n q_{n-1}^{(\lambda)}\|_\infty + \|q_{n-1}^{(\lambda)} x_{n-1} q_{n-1}^{(\lambda)}\|_\infty$$

$$\leq (k+1) \|q_{n-1}^{(\lambda)} x_{n-1} q_{n-1}^{(\lambda)}\|_\infty \leq (k+1)\lambda.$$

Therefore,

$$II \leq 2(k+1)^2 \lambda^2 \tau(\mathbf{1} - q_N^{(\lambda)})$$

$$\leq 2(k+1)^2 \lambda.$$

Hence, $\|q_1^{(\lambda)} dx_1 q_1^{(\lambda)}\|_2^2 + \sum_{n \geq 2} \|q_{n-1}^{(\lambda)} dx_n q_{n-1}^{(\lambda)}\|_2^2 \leq 4\lambda + 2(k+1)^2 \lambda$. Combining all the above estimates, we conclude that

$$\tau(\chi_{(\lambda,\infty)}(\sigma_C(y))) \leq [34 + 16(k+1)^2] \lambda^{-1}.$$

The proof for $\|\sigma_R(z)\|_{1,\infty}$ is identical. $\square$

**4. Best constants for noncommutative Burkholder inequalities.** The following is the principal result of this section.

THEOREM 4.1. *Let $1 < p < \infty$. There exist two constants $\delta_p > 0$ and $\eta_p > 0$ (depending only on $p$) such that for any finite martingale $x$ in $L^p(\mathcal{M}, \tau)$,*

(B$_p$) $$\delta_p^{-1} \|x\|_{h^p(\mathcal{M})} \leq \|x\|_p \leq \eta_p \|x\|_{h^p(\mathcal{M})}.$$

*Moreover, we have the following estimates for the best constants in (B$_p$):*

(i) $\delta_p \approx (p-1)^{-1}$ *as* $p \to 1$;
(ii) $\delta_p \approx p$ *as* $p \to \infty$;
(iii) $\eta_p \approx 1$ *as* $p \to 1$;
(iv) *there exists an absolute constant $C$ such that $\eta_p \leq Cp$ for $p > 2$.*

The inequalities ($B_p$), known as the noncommutative Burkholder inequalities, were originally proved by Junge and Xu in [16]. The main purpose here is to provide the right order of growth for the best constants. The proof in [16] started from establishing the case $p \geq 2$ and then deduced the case $1 < p < 2$ through duality argument. Our approach follows the opposite direction: first we deduce the case $1 < p < 2$ using Theorem 3.1 combined with real interpolations, then deduce the range $p \geq 2$ by duality. This approach has advantages as it leads to the orders of growth of the constants as stated in Theorem 3.1. We treat the two cases $1 < p \leq 2$ and $p > 2$ in two separate subsections.

### 4.1. Noncommutative Burkholder inequalities for $1 < p \leq 2$.

PROPOSITION C. *Let $1 < p < 2$. There exists two constants $\delta_p > 0$ and $\eta_p > 0$ (depending only on $p$) such that for any finite martingale $x$ in $L^p(\mathcal{M}, \tau)$,*

$$(\mathrm{B}_p) \qquad \delta_p^{-1} \|x\|_{h^p(\mathcal{M})} \leq \|x\|_p \leq \eta_p \|x\|_{h^p(\mathcal{M})}.$$

*Moreover, $\delta_p \approx (p-1)^{-1}$ as $p \to 1$ and $\eta_p \approx 1$ as $p \to 1$. These orders of growth are optimal.*

We remark that $\eta_p \approx 1$ as $p \to 1$ was already obtained in [16]. Moreover, it is also noted in [16] that for $1 < p < 2$, $\|x\|_{\mathcal{H}^p} \leq 2^{1/p} \|x\|_{h^p}$ so from Theorem 2.3, we have $\alpha_p \leq 2^{1/p} \delta_p$ and therefore we obtain the estimate of $\delta_p$ from below. Thus, it remains to estimate $\delta_p$ from above. This is a direct application of Theorem 3.1 via interpolation.

LEMMA 4.2. *There exists an absolute constant $C$ such that for $1 < p < 2$, $\delta_p \leq C(p-1)^{-1}$.*

Our main tool is real interpolation, principally the $J$-method. We will review the general setup of the $J$-method. Our main reference for facts about interpolation is the book [1].

A pair of (quasi)-Banach spaces $(E_0, E_1)$ is called a compatible couple if they embed continuously into some topological vector space $X$. This allows us to consider the spaces $E_0 \cap E_1$ and $E_0 + E_1$ equipped with $\|x\|_{E_0 \cap E_1} = \max\{\|x\|_{E_0}, \|x\|_{E_1}\}$, $\|x\|_{E_0 + E_1} = \inf\{\|x_0\|_{E_0} + \|x_1\|_{E_1} : x = x_0 + x_1, x_0 \in E_0, x_1 \in E_1\}$, respectively.

For a compatible couple $(E_0, E_1)$, we define for any $x \in E_0 \cap E_1$, and $t > 0$,

$$J(x, t; E_0, E_1) = \max\{\|x\|_{E_0}, t\|x\|_{E_1}\}.$$

If the compatible couple $(E_0, E_1)$ is clear from the context, we will simply write $J(x, t)$ in place of $J(x, t; E_0, E_1)$.

We will work with the *discrete version of the J-method* which we will now describe: for $0 < \theta < 1$ and $1 \leq p < \infty$, we denote by $\lambda^{\theta,p}$ the space of all sequences $(\alpha_\nu)_{\nu=-\infty}^{\infty}$ for which

$$\|(\alpha_\nu)\|_{\lambda^{\theta,p}} = \left\{\sum_{\nu \in \mathbb{Z}} (2^{-\nu\theta}|\alpha_\nu|)^p\right\}^{1/p} < \infty.$$

DEFINITION 4.3. Let $(E_0, E_1)$ be a compatible couple and suppose that $0 < \theta < 1$, and $1 \leq p < \infty$. The interpolation space $(E_0, E_1)_{\theta,p,J}$ consists of elements $x \in E_0 + E_1$ which admits a representation

(4.1) $$x = \sum_{\nu \in \mathbb{Z}} u_\nu \quad \text{(convergence in } E_0 + E_1\text{)},$$

with $u_\nu \in E_0 \cap E_1$ and such that

$$\|x\|_{\theta,p,J} = \inf\{\|\{J(u_\nu, 2^\nu)\}\|_{\lambda^{\theta,p}}\} < \infty,$$

where the infimum is taken over all representations of $x$ as in (4.1).

PROOF OF LEMMA 4.2. It suffices to verify the lemma for positive finite $L^2$-martingale $x = (x_n)_{1 \leq n \leq N}$. For $1 < p < 2$, let $0 < \theta < 1$ such that $1/p = (1-\theta) + \theta/2$. For $\varepsilon > 0$, fix $(u_\nu)_{\nu=-\infty}^{\infty}$ in $L^2(\mathcal{M}, \tau)$ such that

$$x_N = \sum_{\nu \in \mathbb{Z}} u_\nu$$

and

$$\|x_N\|_{\theta,p;J} + \varepsilon \geq \|\{J(u_\nu, 2^\nu)\}\|_{\lambda^{\theta,p}},$$

where the $J$-functional and the interpolation are relative to the interpolation couple $(L^1(\mathcal{M}), L^2(\mathcal{M}))$.

For each $\nu \in \mathbb{Z}$, Theorem 3.1 guarantees the existence of an absolute constant $K > 0$, and three finite adapted sequences $a^{(\nu)}$, $b^{(\nu)}$ and $c^{(\nu)}$ in $L^2(\mathcal{M}, \tau)$ such that:

(1) $\mathcal{E}_n(u_\nu) - \mathcal{E}_{n-1}(u_\nu) = a_n^{(\nu)} + b_n^{(\nu)} + c_n^{(\nu)}$ for all $1 \leq n \leq N$;
(2) $J(\sum_{n \geq 1} a_n^{(\nu)} \otimes e_{n,n}, t) \leq K J(u_\nu, t)$ for every $t > 0$;
(3) $J((\sum_{n \geq 1} \mathcal{E}_{n-1}(|b_n^{(\nu)}|^2))^{1/2}, t) \leq K J(u_\nu, t)$ for every $t > 0$;
(4) $J((\sum_{n \geq 1} \mathcal{E}_{n-1}(|c_n^{(\nu)*}|^2))^{1/2}, t) \leq K J(u_\nu, t)$ for every $t > 0$,

where the $J$-functionals in the left-hand side of the inequality in (2) above are taken relative to the interpolation couple $(L^{1,\infty}(\mathcal{M}\overline{\otimes}B(l^2)), L^2(\mathcal{M}\overline{\otimes}B(l^2)))$

and those from the left hand sides of (3) and (4) are taken with respect to the interpolation couple $(L^{1,\infty}(\mathcal{M}), L^2(\mathcal{M}))$. From this, we can deduce that

$$\text{(4.2)} \qquad \left\|\left\{J\left(\sum_{n\geq 1} a_n^{(\nu)} \otimes e_{n,n}, 2^\nu\right)\right\}\right\|_{\lambda^{\theta,p}} \leq K(\|x_N\|_{\theta,p;\underline{J}} + \varepsilon),$$

$$\text{(4.3)} \qquad \left\|\left\{J\left(\left(\sum_{n\geq 1} \mathcal{E}_{n-1}(|b_n^{(\nu)}|^2)\right)^{1/2}, 2^\nu\right)\right\}\right\|_{\lambda^{\theta,p}} \leq K(\|x_N\|_{\theta,p;\underline{J}} + \varepsilon)$$

and

$$\text{(4.4)} \qquad \left\|\left\{J\left(\left(\sum_{n\geq 1} \mathcal{E}_{n-1}(|c_n^{(\nu)*}|^2)\right)^{1/2}, 2^\nu\right)\right\}\right\|_{\lambda^{\theta,p}} \leq K(\|x_N\|_{\theta,p;\underline{J}} + \varepsilon).$$

From (4.2) and the definition of $\|\cdot\|_{\theta,p;\underline{J}}$, we get that for any finite subset $S \subset \mathbb{Z}$,

$$\text{(4.5)} \qquad \left\|\sum_{\nu \in S}\sum_{n\geq 1} a_n^{(\nu)} \otimes e_{n,n}\right\|_{[L^{1,\infty}(\mathcal{M}\overline{\otimes}B(l^2)), L^2(\mathcal{M}\overline{\otimes}B(l^2))]_{\theta,p;\underline{J}}} \leq K(\|x_N\|_{\theta,p;\underline{J}} + \varepsilon),$$

and therefore the series $\sum_{\nu\in\mathbb{Z}}(\sum_{n\geq 1} a_n^{(\nu)} \otimes e_{n,n})$ is (unconditionally) convergent in the Banach space $[L^{1,\infty}(\mathcal{M}\overline{\otimes}B(l^2)), L^2(\mathcal{M}\overline{\otimes}B(l^2))]_{\theta,p;\underline{J}}$; hence if we set

$$a := \sum_{\nu\in\mathbb{Z}} a^{(\nu)},$$

then the sequence $a = (a_n)_{n\geq 1}$ satisfies

$$\text{(4.6)} \qquad \left\|\sum_{n\geq 1} a_n \otimes e_{n,n}\right\|_{[L^{1,\infty}(\mathcal{M}\overline{\otimes}B(l^2)), L^2(\mathcal{M}\overline{\otimes}B(l^2))]_{\theta,p;\underline{J}}} \leq K(\|x_N\|_{\theta,p;\underline{J}} + \varepsilon).$$

For the other sequences, we remark from the definition of the $J$-functionals that

$$J\left(\left(\sum_{n\geq 1} \mathcal{E}_{n-1}(|b_n^{(\nu)}|^2)\right)^{1/2}, 2^\nu\right)$$

$$= \max\left\{\left\|\left(\sum_{n\geq 1} \mathcal{E}_{n-1}(|b_n^{(\nu)}|^2)\right)^{1/2}\right\|_{1,\infty}, 2^\nu \left\|\left(\sum_{n\geq 1} \mathcal{E}_{n-1}(|b_n^{(\nu)}|^2)\right)^{1/2}\right\|_2\right\}.$$

Now we consider the conditioned space $L^2_{\text{cond}}(\mathcal{M}; l^2_C)$ as a subspace of the space $L^2(\mathcal{M}\overline{\otimes}B(l^2(\mathbb{N}^2)))$ according to [15] and view $(b^{(\nu)})$ as a column vec-tor with entries from $L^2(\mathcal{M}\overline{\otimes}B(l^2(\mathbb{N}^2)))$. Then for every $\nu \in \mathbb{Z}$,

$$J\left(\left(\sum_{n\geq 1}\mathcal{E}_{n-1}(|b_n^{(\nu)}|^2)\right)^{1/2},2^\nu\right)$$
$$=J(b^{(\nu)},2^\nu;L^{1,\infty}(\mathcal{M}\overline{\otimes}B(l^2(\mathbb{N}^2))),L^2(\mathcal{M}\overline{\otimes}B(l^2(\mathbb{N}^2)))).$$

Then (4.3) becomes

$$(4.7) \quad \begin{aligned}&\|\{J(b^{(\nu)},2^\nu;L^{1,\infty}(\mathcal{M}\overline{\otimes}B(l^2(\mathbb{N}^2))),L^2(\mathcal{M}\overline{\otimes}B(l^2(\mathbb{N}^2))))\}\|_{\lambda^{\theta,p}}\\ &\leq K(\|x_N\|_{\theta,p;\underline{J}}+\varepsilon).\end{aligned}$$

Set $b:=\sum_{\nu\in\mathbb{Z}}b^{(\nu)}$. Then $b$ is a double indexed sequence $(b_{n,k})$ with $b_{n,k}\in L^1(\mathcal{M}_n,\tau|_{\mathcal{M}_n})$ for all $k\in\mathbb{N}$ and satisfies

$$(4.8) \quad \|b\|_{[L^{1,\infty}(\mathcal{M}\overline{\otimes}B(l^2(\mathbb{N}^2))),L^2(\mathcal{M}\overline{\otimes}B(l^2(\mathbb{N}^2)))]_{\theta,p;\underline{J}}}\leq K(\|x_N\|_{\theta,p;\underline{J}}+\varepsilon).$$

We note that a similar argument can be applied to the finite sequences $c^{(\nu)}$'s. That is, if $c:=\sum_{\nu\in\mathbb{Z}}c^{(\nu)}$, then as a row vector, we have

$$(4.9) \quad \|c\|_{[L^{1,\infty}(\mathcal{M}\overline{\otimes}B(l^2(\mathbb{N}^2))),L^2(\mathcal{M}\overline{\otimes}B(l^2(\mathbb{N}^2)))]_{\theta,p;\underline{J}}}\leq K(\|x_N\|_{\theta,p;\underline{J}}+\varepsilon).$$

We remark that since for every $\nu \in \mathbb{Z}$, the sequences $a^{(\nu)}$, $b^{(\nu)}$ and $c^{(\nu)}$ are adapted, it follows that $a$, $b$ and $c$ are adapted sequences. Moreover, it is clear from the construction that for $1\leq n \leq N$,

$$dx_n = a_n + b_n + c_n.$$

We proceed by invoking a general fact about interpolations of noncommutative spaces. First, we recall from the general equivalence theorem on real interpolations that the same inequalities as in (4.6), (4.8), (4.9) can be made with any real interpolation method (with possible change on the absolute constant). Second, it is now understood that for any semifinite von Neumann algebra $\mathcal{N}$ equipped with a semifinite normal trace $\varphi$, the following interpolation results hold:

$$[L^{1,\infty}(\mathcal{N},\varphi),L^2(\mathcal{N},\varphi)]_{\theta,p}=L^p(\mathcal{N},\varphi) \qquad \text{(with equivalent norms)}$$

and

$$[L^1(\mathcal{N},\varphi),L^2(\mathcal{N},\varphi)]_{\theta,p}=L^p(\mathcal{N},\varphi) \qquad \text{(with equivalent norms)}.$$

More precisely, [28], Corollary 2.2, page 1467 implies that it is enough to track the order of the constants for the commutative case. In order to achieve this, we need a few facts. For $f\in L^2$,

$$(4.10) \qquad C(1-\theta)^{-1/2}\|f\|_{L^p}\leq \|f\|_{[L^{1,\infty},L^2]_{\theta,p,K}}$$

and

(4.11) $$\|f\|_{[L^{1,\infty},L^2]_{\theta,p,K}} \leq c(\theta)\|f\|_{[L^{1,\infty},L^2]_{\theta,p,J}},$$

where $c(\theta) = \int_0^\infty s^\theta \min(1, s^{-1})\,ds/s = \theta^{-1}(1-\theta)^{-1}$. The first inequality can be deduced from [13], Theorem 4.3, while the second is in [1], pages 44–45. On the other hand, it is implicit in the proof of [1], Theorem 5.2.1, pages 109–110 that if $1/p + 1/q = 1$, then for every $g \in L^\infty$,

$$\|g\|_q \leq 2\|g\|_{[L^\infty,L^2]_{1-\theta,q,K}}.$$

By duality, we have that for every $f \in L^2$, $\|f\|_{[L^\infty,L^2]^*_{1-\theta,q,K}} \leq 2\|f\|_p$. We now appeal to general duality between the $K$-method and the $J$-method to conclude that for every $f \in L^2$,

(4.12) $$\|f\|_{\theta,p,J} \leq \|f\|_{[L^\infty,L^2]^*_{1-\theta,q,K}} \leq 2\|f\|_p$$

(we note that the constant 1 in the first inequality follows from dualizing the first part of the proof of [1], Theorem 3.7.1, pages 54–55).

Combining (4.6) and (4.8)–(4.12) we can conclude that there exists an absolute constant $C > 0$ such that

$$\left\|\sum_{n\geq 1} a_n \otimes e_{n,n}\right\|_{L^p(\mathcal{M}\overline{\otimes}B(l^2))} + \|b\|_{L^p(\mathcal{M}\overline{\otimes}B(l^2(\mathbb{N}^2)))} + \|c\|_{L^p(\mathcal{M}\overline{\otimes}B(l^2(\mathbb{N}^2)))}$$
$$\leq C\theta^{-1}(\|x\|_p + \varepsilon).$$

From the construction of $a$, $b$ and $c$, this is equivalent to the existence of an absolute constant $C > 0$ such that

(4.13) $$\left(\sum_{n\geq 1}\|a_n\|_p^p\right)^{1/p} + \left\|\left(\sum_{n\geq 1}\mathcal{E}_{n-1}(|b_n^2|)\right)^{1/2}\right\|_p + \left\|\left(\sum_{n\geq 1}\mathcal{E}_{n-1}(|c_n^{*\,2}|)\right)^{1/2}\right\|_p$$
$$\leq C(p-1)^{-1}(\|x\|_p + \varepsilon).$$

We remark that the sequences $a$, $b$ and $c$ are adapted but are not necessarily martingale difference sequences. To complete the proof, it is enough to set for $n \geq 1$,

$$dx_n^D = a_n - \mathcal{E}_{n-1}(a_n),$$
$$dx_n^C = b_n - \mathcal{E}_{n-1}(b_n),$$
$$dx_n^R = c_n - \mathcal{E}_{n-1}(c_n).$$

Then $(dx_n^D)_{n\geq 1}$, $(dx_n^C)_{n\geq 1}$ and $(dx_n^C)_{n\geq 1}$ are martingale difference sequences with $dx_n = dx_n^D + dx_n^C + dx_n^R$ for $1 \leq n \leq N$. The fact that any conditional expectation $\mathcal{E}$ is a contractive projection in $L^p(\mathcal{M}, \tau)$ and satisfies

$\mathcal{E}(y)^*\mathcal{E}(y) \leq \mathcal{E}(y^*y)$ implies that from (4.13) we can deduce

$$\left(\sum_{n\geq 1}\|dx_n^D\|_p^p\right)^{1/p} + \left\|\left(\sum_{n\geq 1}\mathcal{E}_{n-1}(|dx_n^C|^2)\right)^{1/2}\right\|_p$$

(4.14)
$$+ \left\|\left(\sum_{n\geq 1}\mathcal{E}_{n-1}(|dx_n^{R^*}|^2)\right)^{1/2}\right\|_p$$

$$\leq C'(p-1)^{-1}(\|x\|_p + \varepsilon)$$

with $C' = 2C$. Taking the infimum over $\varepsilon > 0$, we conclude from the definition of the $h^p$-norm that

$$\|x\|_{h^p(\mathcal{M})} \leq C'(p-1)^{-1}\|x\|_p,$$

which shows that $\delta_p \leq C'(p-1)^{-1}$ for $1 < p < 2$. Thus the proof is complete. □

4.2. *Noncommutative Burkholder inequalities for $p \geq 2$.*

PROPOSITION D. *Let $2 \leq p < \infty$. There exist two constants $\delta_p > 0$ and $\eta_p > 0$ (depending only on $p$) such that for any finite martingales in $L^p(\mathcal{M}, \tau)$,*

(4.15) $$\delta_p^{-1}\|x\|_{h^p(\mathcal{M})} \leq \|x\|_p \leq \eta_p\|x\|_{h^p(\mathcal{M})}.$$

*Moreover, $\delta_p \approx p$ as $p \to \infty$ and $\eta_p \leq Cp$ as $p \to \infty$ (for some absolute constant $C$).*

The fact that $\delta_p \geq Cp$ is already known from [19], Remark 10. For the remaining estimates, we will use duality arguments. The main ingredient here is a conditioned version of the duality between column (resp. row) spaces developed in [16] which we now state.

LEMMA 4.4 ([16], Lemma 6.5). *Let $1 < p < \infty$ and $1/p + 1/p' = 1$. Then for any $b \in L^{p'}_{\text{cond}}(\mathcal{M}; l_C^2)$, the functional $\xi_b : L^p_{\text{cond}}(\mathcal{M}; l_C^2) \to \mathbb{C}$ defined by $\xi_b(a) = \sum_{n\geq 1}\tau(b_n^* a_n)$ is continuous with*

$$\|\xi_b\| \leq \|b\|_{L^{p'}_{\text{cond}}(\mathcal{M};l_C^2)} \leq \gamma_{p'}\|\xi_b\|,$$

*where $\gamma_{p'}$ is the constant from the noncommutative Stein inequality [27]. Conversely, any functional $\xi \in (L^p_{\text{cond}}(\mathcal{M}; l_C^2))^*$ is given by some sequence $b$ in $L^{p'}_{\text{cond}}(\mathcal{M}; l_C^2)$. A similar statement holds for conditioned row spaces.*

LEMMA 4.5 ([16], Lemma 6.4). *Let $1 \leq p < \infty$. For any finite sequence $a = (a_n)_{n \geq 1}$ in $\mathcal{M}$, define*

$$R(a) = (\mathcal{E}_n(a_n))_{n \geq 1} \quad \text{and} \quad R'(a) = (\mathcal{E}_{n-1}(a_n))_{n \geq 1}.$$

*Then $R$ and $R'$ extend to contractive projections on $L^p_{\text{cond}}(\mathcal{M}; l^2_C)$ and $L^p_{\text{cond}}(\mathcal{M}; l^2_R)$. Consequently, $h^p_C(\mathcal{M})$ [resp. $h^p_R(\mathcal{M})$] is a 2-complemented subspace of $L^p_{\text{cond}}(\mathcal{M}; l^2_C)$ [resp. $L^p_{\text{cond}}(\mathcal{M}; l^2_R)$].*

Using Lemmas 4.4 and 4.5, it is straightforward to verify that if $2 \leq p < \infty$ and $1/p + 1/p' = 1$, then $\delta_p \leq 2\gamma_p \eta_{p'}$ and $\eta_p \leq \delta_{p'}$. The conclusion follows from the facts that $\gamma_p \approx p$ as $p \to \infty$ (see [19]) and $\eta_{p'} \approx 1$ and $\delta_{p'} \approx p$ as $p \to \infty$ from Proposition C.

We conclude this section with the following remark.

REMARK 4.6. (i) Denote by $\eta_p^{(\text{com})}$ and $\delta_p^{(\text{com})}$ the corresponding best constants for the commutative case. We recall that when $p \to \infty$, the optimal order of growth of $\eta_p^{(\text{com})}$ and $\delta_p^{(\text{com})}$ is $0(p/\log p)$ and $0(\sqrt{p})$, respectively. The first is a result of Johnson, Schechtman and Zinn from [14] for the case of $p$-moments of sums of independent mean-zero random variables which was generalized by Hitczenko [12] for the more general case of martingales. The second also appeared in [12]. Thus at the time of this writing, the exact order of growth of the constant $\eta_p$ (when $p \to \infty$) is still open. We can only state the existence of absolute constants $C_1$ and $C_2$ such that $C_1 p/(\log p) \leq \eta_p \leq C_2 p$ when $p$ is large enough.

(ii) If we denote by $\eta_p^{(\text{Ros})}$ and $\delta_p^{(\text{Ros})}$ the corresponding best constants for the case of sums of (noncommutative) independent sequences, then the estimates $\eta_p^{(\text{Ros})} \leq Cp$ and $\delta_p^{(\text{Ros})} \leq Cp$ (for some absolute constant $C$) when $p > 2$ were also obtained in [16].

**5. Noncommutative Rosenthal inequalities and BMO-spaces.** We start by recalling the definitions of BMO-spaces for noncommutative martingales introduced in [27]. Let

$$\text{BMO}_C(\mathcal{M}) := \left\{ a \in L^2(\mathcal{M}, \tau) : \sup_{n \geq 1} \|\mathcal{E}_n|a - \mathcal{E}_{n-1}a|^2\|_\infty < \infty \right\}.$$

Then $\text{BMO}_C(\mathcal{M})$ becomes a Banach space when equipped with the norm

$$\|a\|_{\text{BMO}_C} = \left( \sup_{n \geq 1} \|\mathcal{E}_n|a - \mathcal{E}_{n-1}a|^2\|_\infty \right)^{1/2}.$$

Similarly, we define $\text{BMO}_R(\mathcal{M})$ as the space of all $a$ with $a^* \in \text{BMO}_C(\mathcal{M})$ equipped with the natural norm $\|a\|_{\text{BMO}_R(\mathcal{M})} = \|a^*\|_{\text{BMO}_C(\mathcal{M})}$. The space $\text{BMO}(\mathcal{M})$ is the intersection of these two spaces:

$$\text{BMO}(\mathcal{M}) := \text{BMO}_C(\mathcal{M}) \cap \text{BMO}_R(\mathcal{M})$$

with the intersection norm

$$\|a\|_{\mathrm{BMO}} = \max\{\|a\|_{\mathrm{BMO}_C}, \|a\|_{\mathrm{BMO}_R}\}.$$

We recall that as in the classical case, for $1 < p < \infty$, $\mathcal{M} \subset \mathrm{BMO}(\mathcal{M}) \subset L^p(\mathcal{M}, \tau)$. For more information on martingale BMO-spaces, we refer to [16, 17, 23, 27].

In this section, we consider the case of independence studied in [18, 35].

DEFINITION 5.1. (i) Let $\mathcal{N}$ be a von Neumann subalgebra of $\mathcal{M}$ and let $\mathcal{E}_{\mathcal{N}}$ be the associated trace-preserving normal conditional expectation onto $\mathcal{N}$. A sequence $(\mathcal{A}_n)_{n \geq 1}$ of von Neumann subalgebras of $\mathcal{M}$ is called order independent with respect to $\mathcal{E}_{\mathcal{N}}$ (or $\mathcal{N}$) if for every $n \geq 1$, every $a \in \mathcal{A}_n$ and $b$ in the von Neumann subalgebra generated by $(\mathcal{A}_1, \ldots, \mathcal{A}_{n-1})$,

$$\mathcal{E}_{\mathcal{N}}(ab) = \mathcal{E}_{\mathcal{N}}(a)\mathcal{E}_{\mathcal{N}}(b).$$

(ii) A sequence $(a_n)_{n \geq 1}$ in $L^p(\mathcal{M}, \tau)$ ($2 \leq p \leq \infty$) is called independent with respect to $\mathcal{E}_{\mathcal{N}}$ if there is an order-independent sequence $(\mathcal{A}_n)_{n \geq 1}$ of von Neumann subalgebras of $\mathcal{M}$ such that $a_n \in L^p(\mathcal{A}_n)$ for all $n \geq 1$.

If $\mathcal{N} = \mathbb{C}\mathbf{1}$ [then $\mathcal{E}_{\mathcal{N}} = \tau(\cdot)\mathbf{1}$], we simply say "independent" with respect to $\tau$.

We refer to [18, 35] for natural examples of independent sequences.

REMARK 5.2. Let $(\mathcal{A}_n)_{n \geq 1}$ be an independent sequence of von Neumann subalgebras and $a_n \in L^p(\mathcal{A}_n)$ with $\mathcal{E}_{\mathcal{N}}(a_n) = 0$. Then $(a_n)_{n \geq 1}$ is a martingale difference sequence in $L^p(\mathcal{M}, \tau)$.

Indeed, if for $n \geq 1$, we set $\mathcal{M}_n$ to be the von Neumann subalgebra generated by $(\mathcal{A}_1, \ldots, \mathcal{A}_n)$, then $(\mathcal{M}_n)_{n \geq 1}$ is an increasing filtration of von Neumann subalgebras of $\mathcal{M}$. Let $\mathcal{E}_n$ be the associated conditional expectations. Then the independence assumption implies that for every $b \in \mathcal{M}_{n-1}$, $\mathcal{E}_{\mathcal{N}}(\mathcal{E}_{n-1}(a_n)b) = \mathcal{E}_{\mathcal{N}}(a_n b) = \mathcal{E}_{\mathcal{N}}(a_n)\mathcal{E}_{\mathcal{N}}(b)$. Therefore, $\mathcal{E}_{n-1}(a_n) = \mathcal{E}_{\mathcal{N}}(a_n) = 0$. Thus $(a_n)_{n \geq 1}$ is a martingale difference sequence with respect to the filtration $(\mathcal{M}_n)_{n \geq 1}$.

The next result can be viewed as an extension of the noncommutative Rosenthal inequalities from [18, 35] to the case $p = \infty$. A precursor of this result for commuting sequences can be found in [22].

THEOREM 5.3. *Let $\mathcal{N}$ be a von Neumann subalgebra of $\mathcal{M}$ with the associated normal conditional expectation $\mathcal{E}_{\mathcal{N}}$. Let $(a_n)_{n \geq 1} \subset \mathcal{M}$ be an independent sequence with respect to $\mathcal{E}_{\mathcal{N}}$ such that $\mathcal{E}_{\mathcal{N}}(a_n) = 0$. Then*

$$\left\|\sum_{n \geq 1} a_n\right\|_{\mathrm{BMO}} \sim_C \sup_{n \geq 1}\|a_n\|_{\infty} + \left\|\left(\sum_{n \geq 1} \mathcal{E}_{\mathcal{N}}(a_n a_n^* + a_n a_n^*)\right)^{1/2}\right\|_{\infty}.$$

*In particular, if $\mathcal{N} = \mathbb{C}\mathbf{1}$, then*

$$\left\|\sum_{n\geq 1} a_n\right\|_{\mathrm{BMO}} \sim_C \sup_{n\geq 1} \|a_n\|_\infty + \left(\sum_{n\geq 1} \|a_n\|_2^2\right)^{1/2}.$$

*(Here the BMO-norm is relative to the filtration described above.)*

PROOF. First, we consider a general inequality for martingales. If $x = (x_n)_{n\geq 1}$ is a noncommutative martingale in $L^2(\mathcal{M}, \tau)$, then (see, e.g., [27]) for every $n \geq 1$,

$$\mathcal{E}_n|x_\infty - x_{n-1}|^2 = \mathcal{E}_n\left(\sum_{k\geq n} |dx_k|^2\right)$$

(5.1)

$$= |dx_n|^2 + \mathcal{E}_n\left(\sum_{k\geq n+1} \mathcal{E}_{k-1}(|dx_k|^2)\right).$$

Taking the norm in $\mathcal{M}$ and supremum over $n \geq 1$, we deduce that

$$\|x\|_{\mathrm{BMO}_C} \leq \sup_{n\geq 1} \|dx_n\|_\infty + \|\sigma_C(x)\|_\infty.$$

Observe that $\sigma_C(x)^2 \leq \sum_{n\geq 1} \mathcal{E}_{n-1}(|dx_n|^2 + |dx_n^*|^2)$. Combining with similar argument for the $\mathrm{BMO}_R$-norm, we have

$$(5.2) \quad \|x\|_{\mathrm{BMO}} \leq \sup_{n\geq 1} \|dx_n\|_\infty + \left\|\left(\sum_{n\geq 1} \mathcal{E}_{n-1}(|dx_n|^2 + |dx_n^*|^2)\right)^{1/2}\right\|_\infty.$$

Now, let $(a_n)_{n\geq 1}$ be an independent sequence with respect to $\mathcal{E}_\mathcal{N}$. Then by the remark above $(a_n)_{n\geq 1}$ is a martingale difference sequence and therefore (5.2) applies. Moreover, the independence assumption implies that for every $n \geq 1$,

$$\mathcal{E}_{n-1}(a_n^* a_n + a_n^* a_n) = \mathcal{E}_\mathcal{N}(a_n^* a_n + a_n^* a_n).$$

Hence (5.2) becomes

$$\left\|\sum_{n\geq 1} a_n\right\|_{\mathrm{BMO}} \leq \sup_{n\geq 1} \|a_n\|_\infty + \left\|\left(\sum_{n\geq 1} \mathcal{E}_\mathcal{N}(a_n^* a_n + a_n a_n^*)\right)^{1/2}\right\|_\infty.$$

Thus one inequality is proved. For the reverse inequality, we note that if $a = \sum_{n\geq 1} a_n$, then from (5.1) we have

$$\mathcal{E}_n|a - \mathcal{E}_{n-1}a|^2 = |a_n|^2 + \mathcal{E}_\mathcal{N}\left(\sum_{k\geq n+1} a_k^* a_k\right).$$

NONCOMMUTATIVE MARTINGALES 29

In particular, $\mathcal{E}_n|a - \mathcal{E}_{n-1}a|^2 \geq |a_n|^2$ and therefore
$$\|a\|_{\mathrm{BMO}_C} \geq \sup_n \|a_n\|_\infty.$$

Moreover,
$$\mathcal{E}_\mathcal{N}\left(\sum_{n\geq 1} a_n^* a_n\right) = \mathcal{E}_\mathcal{N}(a_1^* a_1) + \sum_{n\geq 2} \mathcal{E}_\mathcal{N}(a_n^* a_n)$$
$$\leq \mathcal{E}_\mathcal{N}(a_1^* a_1) + \mathcal{E}_1|a - \mathcal{E}_0 a|^2.$$

Hence,
$$\left\|\sum_{n\geq 1} \mathcal{E}_\mathcal{N}(a_n^* a_n)\right\|_\infty \leq \|a_1\|_\infty^2 + \|a\|_{\mathrm{BMO}_C}^2.$$

Similar argument with the adjoint operators gives
$$\left\|\sum_{n\geq 1} \mathcal{E}_\mathcal{N}(a_n a_n^*)\right\|_\infty \leq \|a_1\|_\infty^2 + \|a\|_{\mathrm{BMO}_R}^2.$$

Combining the last two inequalities, we have
$$\left\|\sum_{n\geq 1} \mathcal{E}_\mathcal{N}(a_n^* a_n + a_n a_n^*)\right\|_\infty \leq 2\|a_1\|_\infty^2 + 2\|a\|_{\mathrm{BMO}}^2.$$

We can now conclude that
$$\sup_{n\geq 1} \|a_n\|_\infty + \left\|\left(\sum_{n\geq 1} \mathcal{E}_\mathcal{N}(a_n^* a_n + a_n a_n^*)\right)^{1/2}\right\|_\infty \leq (1+\sqrt{3})\|a\|_{\mathrm{BMO}}.$$

Thus the proof of the theorem is complete. □

Following [35], from the preceding theorem, we can deduce the following Khintchine-type inequality relative to BMO-spaces.

COROLLARY 5.4. *Let $(a_n)_{n\geq 1}$ be an independent (with respect to $\tau$) sequence in $\mathcal{M}$ with $\tau(a_n) = 0$ for all $n \geq 1$. Assume that*
$$\inf_{n\geq 1} \|a_n\|_2 = \alpha > 0 \quad and \quad \sup_{n\geq 1} \|a_n\|_\infty = \beta < \infty.$$

*Let $\mathcal{B}$ be a finite von Neumann algebra equipped with normal tracial state $\nu$. Then for any finite sequence $b = (b_n)_{n\geq 1}$ in $\mathcal{B}$,*
$$\left\|\sum_{n\geq 1} a_n \otimes b_n\right\|_{\mathrm{BMO}(\mathcal{M}\overline{\otimes}\mathcal{B})} \sim_C \|b\|_{L^\infty(\mathcal{B}, l_C^2) \cap L^\infty(\mathcal{B}, l_R^2)}.$$

The proof follows verbatim the argument of [35], Corollary 6.3, and is left to the interested reader. We should compare Corollary 5.4 with Voiculescu's inequality from [34]. This translates into the following statement: let $\mathcal{A} = A_1 * A_2 * \cdots * A_n$ denote the reduced free product of a finite sequence $A_1, A_2, \ldots, A_n$ of von Neumann algebras respectively equipped with tracial states $\tau_1, \tau_2, \ldots, \tau_n$. Then given $a_1 \in A_1$, $a_2 \in A_2, \ldots, a_n \in A_n$ mean-zero (freely independent) random variables in $\mathcal{A}$ and $b_1, b_2, \ldots, b_n \in \mathcal{B}$, then combining Voiculescu's inequality and Corollary 5.4, we have

$$\left\| \sum_{k \geq 1} a_k \otimes b_k \right\|_{\mathrm{BMO}(\mathcal{A} \overline{\otimes} \mathcal{B})} \sim_C \left\| \sum_{k \geq 1} a_k \otimes b_k \right\|_{\mathcal{A} \overline{\otimes} \mathcal{B}}.$$

**6. Concluding remarks and related open problems.** (i) The optimal orders of growth of the constants in Theorem 4.1 remain valid for the more general martingales on Haagerup's $L^p$-spaces using Haagerup's approximation [11]. This follows from a general deduction of martingale inequalities in type-III cases from finite cases, achieved by Junge and Xu (still unpublished notes).

(ii) In [25], another weak-type inequality was considered for classical martingales as extension of Burkholder inequality which we state explicitly:

THEOREM 6.1 ([25], Theorem A). *Let $(\Omega, \mathcal{F}, \mathbb{P})$ be a probability space and let $(f_n)_{n \geq 1}$ be a martingale in $L^1(\Omega, \mathcal{F}, \mathbb{P})$; then there exist two martingales $(g_n)_{n \geq 1}$ and $(h_n)_{n \geq 1}$ with $f_n = g_n + h_n$ for all $n \geq 1$ that satisfy*

$$\left\| \sum_{n \geq 1} |dh_n| \right\|_{1,\infty} + \left\| \left( \sum_{n \geq 1} \mathbb{E}_{n-1}(|dg_n|^2) \right)^{1/2} \right\|_{1,\infty} \leq K \sup_{n \geq 1} \|f_n\|_1.$$

An example given in [25] shows that this formulation is not comparable to the one considered in this paper. The decomposition in Theorem 6.1 is exactly the classical Davis decomposition. It is noted in [25] that this weak-type inequality does not imply the Burkholder inequality through interpolation. In fact from the well-known property of the Davis decomposition, namely $\|\sum_{n \geq 1} |dh_n|\|_p \leq c_p \|f\|_p$ for $1 < p < \infty$ [4], it is probably more accurate to describe Theorem 6.1 as a weak-type $(1,1)$ extension of the classical Davis theorem. It is still unknown if Theorem 6.1 has noncommutative analogues.

(iii) Combining the noncommutative Burkholder–Gundy inequalities (Theorem 2.3) and the noncommutative Burkholder inequalities (Theorem 4.1), we can state:

PROPOSITION 6.2. *Let $1 < p < \infty$. There exist two constants $\kappa_p > 0$ and $\upsilon_p > 0$ (depending only on $p$) such that for any finite martingale $x$* in $L^p(\mathcal{M}, \tau)$,

$$\kappa_p^{-1} \|x\|_{h^p(\mathcal{M})} \leq \|x\|_{\mathcal{H}^p(\mathcal{M})} \leq \upsilon_p \|x\|_{h^p(\mathcal{M})}.$$

*Moreover:*

(i) $\kappa_p = 0((p-1)^{-1})$ *as* $p \to 1$;
(ii) $\kappa_p = 0(p)$ *as* $p \to \infty$;
(iii) $\upsilon_p \approx 1$ *as* $p \to 1$;
(iv) $\upsilon_p = 0(\sqrt{p})$ *as* $p \to \infty$.

We do not know if these orders of growth are optimal.

(iv) We end the paper with a note on the conditioned Hardy space $h^1$. Let us recall the classical Davis theorem for commutative martingales (see [4]):

$$\left\| \left( \sum_n |df_n|^2 \right)^{1/2} \right\|_1 \sim_C \left\| \sup_n |f_n| \right\|_1.$$

If we denote by $\mathcal{H}^1_{\max}$ the space of commutative martingales $(f_n)_{n \geq 1}$ with $\sup_n |f_n| \in L^1$, then the Davis theorem means that in the commutative case

$$\mathcal{H}^1 = \mathcal{H}^1_{\max} \qquad \text{(with equivalent norms).}$$

Note that $h^1 \subset \mathcal{H}^1$. It turns out the (commutative) conditioned Hardy space $h^1$ coincides with the other two Hardy spaces. Indeed, [10], Theorem IV1.2, page 127, together with the classical Davis decomposition imply that there is a constant $C$ such that $\|f\|_{h^1} \leq C \|f\|_{\mathcal{H}^1_{\max}}$ for every martingale $f$. Therefore we can state that for commutative martingales,

$$h^1 = \mathcal{H}^1 = \mathcal{H}^1_{\max} \qquad \text{(with equivalent norms).}$$

The noncommutative case is surprisingly different as noted in [19]. Indeed, it was shown in [19], Corollary 14, that $\mathcal{H}^1$ and $\mathcal{H}^1_{\max}$ do not coincide in general. Motivated by the commutative case, one could ask about the position of the space $h^1$ with respect to $\mathcal{H}^1$ (or $\mathcal{H}^1_{\max}$) for the noncommutative case.

**Acknowledgments.** I would like to thank Gedon Schechtman for bringing the reference [22] to my attention. I am also indebted to the referee for several comments that significantly improve the presentation of the paper.

## REFERENCES

[1] BERGH, J. and LÖFSTRÖM, J. (1976). *Interpolation Spaces. An Introduction.* Springer, Berlin. MR0482275
[2] BURKHOLDER, D. L. (1973). Distribution function inequalities for martingales. *Ann. Probab.* **1** 19–42. MR0365692

[3] CUCULESCU, I. (1971). Martingales on von Neumann algebras. *J. Multivariate Anal.* **1** 17–27. MR0295398

[4] DAVIS, B. (1970). On the integrability of the martingale square function. *Israel J. Math.* **8** 187–190. MR0268966

[5] DIESTEL, J. and UHL JR., J. J. (1977). *Vector Measures.* Amer. Math. Soc., Providence, RI. MR0453964

[6] DIXMIER, J. (1953). Formes linéaires sur un anneau d'opérateurs. *Bull. Soc. Math. France* **81** 9–39. MR0059485

[7] DOOB, J. L. (1953). *Stochastic Processes.* Wiley, New York. MR0058896

[8] EDWARDS, R. E. and GAUDRY, G. I. (1977). *Littlewood–Paley and Multiplier Theory.* Springer, Berlin. MR0618663

[9] FACK, T. and KOSAKI, H. (1986). Generalized $s$-numbers of $\tau$-measurable operators. *Pacific J. Math.* **123** 269–300. MR0840845

[10] GARSIA, A. M. (1973). *Martingale Inequalities*: *Seminar Notes on Recent Progress.* W. A. Benjamin, Inc., Reading, MA. MR0448538

[11] HAAGERUP, U. (1980). Non-commutative integration theory. *Lecture given at the Symposium in Pure Mathematics of the Amer. Math. Soc., Queens University, Kingston, Ontario.*

[12] HITCZENKO, P. (1990). Best constants in martingale version of Rosenthal's inequality. *Ann. Probab.* **18** 1656–1668. MR1071816

[13] HOLMSTEDT, T. (1970). Interpolation of quasi-normed spaces. *Math. Scand.* **26** 177–199. MR0415352

[14] JOHNSON, W. B., SCHECHTMAN, G. and ZINN, J. (1985). Best constants in moment inequalities for linear combinations of independent and exchangeable random variables. *Ann. Probab.* **13** 234–253. MR0770640

[15] JUNGE, M. (2002). Doob's inequality for non-commutative martingales. *J. Reine Angew. Math.* **549** 149–190. MR1916654

[16] JUNGE, M. and XU, Q. (2003). Noncommutative Burkholder/Rosenthal inequalities. *Ann. Probab.* **31** 948–995. MR1964955

[17] JUNGE, M. and MUSAT, M. (2007). A noncommutative version of the John–Nierenberg theorem. *Trans. Amer. Math. Soc.* **359** 115–142. MR2247885

[18] JUNGE, M. and XU, Q. (2005). Noncommutative Burkholder/Rosenthal inequalities II: Applications. Preprint.

[19] JUNGE, M. and XU, Q. (2005). On the best constants in some non-commutative martingale inequalities. *Bull. London Math. Soc.* **37** 243–253. MR2119024

[20] KADISON, K. V. and RINGROSE, J. R. (1983). *Fundamentals of the Theory of Operator Algebras. I.* Academic Press, New York. MR0719020

[21] LINDENSTRAUSS, J. and TZAFRIRI, L. (1979). *Classical Banach Spaces. II.* Springer, Berlin. MR0540367

[22] MÜLLER, P. F. X. and SCHECHTMAN, G. (1989). On complemented subspaces of $H^1$ and VMO. *Geometric Aspects of Functional Analysis. Lecture Notes in Math.* **1376** 113–125. Springer, Berlin. MR1091439

[23] MUSAT, M. (2003). Interpolation between non-commutative BMO and non-commutative $L^p$-spaces. *J. Funct. Anal.* **202** 195–225. MR1994770

[24] NELSON, E. (1974). Notes on non-commutative integration. *J. Funct. Anal.* **15** 103–116. MR0355628

[25] PARCET, J. (2006). Weak type estimates associated to Burkholder's martingale inequality. Preprint.

[26] PARCET, J. and RANDRIANANTOANINA, N. (2006). Gundy's decomposition for non-commutative martingales and applications. *Proc. London Math. Soc.* (*3*) **93** 227–252. MR2235948
[27] PISIER, G. and XU, Q. (1997). Non-commutative martingale inequalities. *Comm. Math. Phys.* **189** 667–698. MR1482934
[28] PISIER, G. and XU, Q. (2003). Non-commutative $L^p$-spaces. In *Handbook of the Geometry of Banach Spaces* **2** 1459–1517. North-Holland, Amsterdam. MR1999201
[29] RANDRIANANTOANINA, N. (2002). Non-commutative martingale transforms. *J. Funct. Anal.* **194** 181–212. MR1929141
[30] RANDRIANANTOANINA, N. (2004). Square function inequalities for non-commutative martingales. *Israel J. Math.* **140** 333–365. MR2054851
[31] RANDRIANANTOANINA, N. (2005). A weak-type inequality for non-commutative martingales and applications. *Proc. London Math. Soc.* (*3*) **91** 509–544. MR2167096
[32] ROSENTHAL, H. P. (1970). On the subspaces of $L^p$ ($p \geq 2$) spanned by sequences of independent random variables. *Israel J. Math.* **8** 273–303. MR0271721
[33] TAKESAKI, M. (1979). *Theory of Operator Algebras. I.* Springer, New York. MR0548728
[34] VOICULESCU, D. (1998). A strengthened asymptotic freeness result for random matrices with applications to free entropy. *Internat. Math. Res. Notices* **1** 41–63. MR1601878
[35] XU, Q. (2003). Recent development on non-commutative martingale inequalities. In *Functional Space Theory and Its Applications. Proceedings of International Conference & 13th Academic Symposium in China* 283–314. Wuhan, Research Information Ltd.

DEPARTMENT OF MATHEMATICS AND STATISTICS
MIAMI UNIVERSITY
OXFORD, OHIO 45056
USA
E-MAIL: randrin@muohio.edu